\pdfoutput=1
\RequirePackage{ifpdf}
\ifpdf 
\documentclass[pdftex]{sigma}
\else
\documentclass{sigma}
\fi

\usepackage{multirow,mdwtab,tikz,empheq,caption,multirow,array,tabularx}

\newcommand{\bea}{\begin{eqnarray}}
\newcommand{\eea}{\end{eqnarray}}
\newcommand{\bee}{\begin{eqnarray*}}
\newcommand{\eee}{\end{eqnarray*}}
\newcommand{\al}{\begin{align*}}
\newcommand{\eal}{\end{align*}}
\newcommand{\be}{\begin{equation}}
\newcommand{\ee}{\end{equation}}
\newcommand{\eq}[1]{(\ref{#1})}
\newcommand{\bem}{\begin{pmatrix}}
\newcommand{\eem}{\end{pmatrix}}

\def\a{\alpha}

\def\d{\delta}

\def\f{\phi}

\def\inf{\infty}

\def\l{\lambda}
\def\m{\mu}

\def\o{\omega}
\def\w{\omega}
\def\p{\pi}
\def\pa{\partial}

\def\s{\sigma}
\def\t{\tau}
\def\th{\theta}
\def\til{\tilde}

\def\z{\zeta}

\def\D{\Delta}

\def\Tr{\operatorname{Tr}}

\def\Gg{H}
\def\Hh{X}
\newcolumntype{R}{ >{$}r <{$}}
\newcolumntype{C}{ >{$}c <{$}}
\newcolumntype{L}{ >{$}l <{$}}
\newcolumntype{F}{>{\centering\arraybackslash}m{1.5cm}}

\def\ll{\ell}

\newcommand{\comment}[1]{}

\newcommand{\RR}{{\mathbb R}}
\newcommand{\CC}{{\mathbb C}}
\newcommand{\ZZ}{{\mathbb Z}}
\newcommand{\QQ}{{\mathbb Q}}
\newcommand{\HH}{{\mathbb H}}

\newcommand{\Sym}{{{\rm Sym}}}

\newcommand{\Dih}{{{\rm Dih}}}

\newcommand{\ex}{\operatorname{e}} 

\newcommand{\PSL}{\operatorname{PSL}} 
\newcommand{\SL}{\operatorname{SL}} 
\newcommand{\PGL}{\operatorname{PGL}} 
\newcommand{\AGL}{{{\rm AGL}}} 
\newcommand{\GL}{{{\rm GL}}} 
\newcommand{\SO}{\operatorname{\SO}} 

\newcommand{\G}{\Gamma}	
\newcommand{\g}{\gamma}	

\newcommand{\rs}{{X}}	

\newcommand{\xmod}{{\rm \;mod\;}}

\newcommand{\ordg}{{\operatorname{ord}(g)}}

\def\Gell{G^{(\ell)}}

\def\bz{\mathbb{Z}}
\newcommand{\sm}[4]{\left( \begin{smallmatrix} #1&#2\\ #3&#4 \end{smallmatrix} \right)}

\numberwithin{equation}{section}

\theoremstyle{definition}

\makeatletter
\newsavebox\myboxA
\newsavebox\myboxB
\newlength\mylenA

\newcommand*\xunderline[2][0.75]{%
 \sbox{\myboxA}{$\m@th#2$}%
 \setbox\myboxB\null
 \ht\myboxB=\ht\myboxA%
 \dp\myboxB=\dp\myboxA%
 \wd\myboxB=#1\wd\myboxA
 \sbox\myboxB{$\m@th\underline{\copy\myboxB}$}
 \setlength\mylenA{\the\wd\myboxA}
 \addtolength\mylenA{-\the\wd\myboxB}%
 \ifdim\wd\myboxB<\wd\myboxA%
 \rlap{\hskip 0.5\mylenA\usebox\myboxB}{\usebox\myboxA}%
 \else
 \hskip -0.5\mylenA\rlap{\usebox\myboxA}{\hskip 0.5\mylenA\usebox\myboxB}%
 \fi}
\makeatother

\newcommand*\widefbox[1]{\fbox{\hspace{1em}#1\hspace{1em}}}

\newcommand\halfbox[1]{
 \tikz[baseline=(n.base)]{\node(n)[inner sep=1pt]{$#1$};
 \draw[line cap=round](n.north west)--(n.south west)--(n.south east);
 }
}

\begin{document}

\allowdisplaybreaks

\newcommand{\arXivNumber}{1608.07835}

\renewcommand{\thefootnote}{}

\renewcommand{\PaperNumber}{014}

\FirstPageHeading

\ShortArticleName{Generalised Umbral Moonshine}

\ArticleName{Generalised Umbral Moonshine\footnote{This paper is a~contribution to the Special Issue on Moonshine and String Theory. The full collection is available at \href{https://www.emis.de/journals/SIGMA/moonshine.html}{https://www.emis.de/journals/SIGMA/moonshine.html}}}

\Author{Miranda C.N.~CHENG~$^{\dag^1\dag^2}$, Paul DE LANGE~$^{\dag^3}$ and Daniel P.Z.~WHALEN~$^{\dag^4}$}

\AuthorNameForHeading{M.C.N.~Cheng, P.~de Lange and D.P.Z.~Whalen}

\Address{$^{\dag^1}$~Korteweg-de Vries Institute for Mathematics, Amsterdam, The Netherlands}
\EmailDD{\href{mailto:mcheng@uva.nl}{mcheng@uva.nl}}

\Address{$^{\dag^2}$~Institute of Physics, University of Amsterdam, Amsterdam, The Netherlands}

\Address{$^{\dag^3}$~Department of Physics and Astronomy, University of Kentucky, Lexington, KY 40506, USA}
\EmailDD{\href{mailto:p.delange@uky.edu}{p.delange@uky.edu}}

\Address{$^{\dag^4}$~Stanford Institute for Theoretical Physics, Department of Physics and Theory Group, SLAC, \\
\hphantom{$^{\dag^4}$}~Stanford University, Stanford, CA 94305, USA}
\EmailDD{\href{mailto:dpzwhalen@gmail.com}{dpzwhalen@gmail.com}}

\ArticleDates{Received October 08, 2018, in final form January 30, 2019; Published online March 02, 2019}

\Abstract{Umbral moonshine describes an unexpected relation between 23 finite groups arising from lattice symmetries and special mock modular forms. It includes the Mathieu moonshine as a special case and can itself be viewed as an example of the more general moonshine phenomenon which connects finite groups and distinguished modular objects. In this paper we introduce the notion of generalised umbral moonshine, which includes the generalised Mathieu moonshine~[Gaberdiel~M.R., Persson~D., Ronellenfitsch~H., Volpato~R., \textit{Commun. Number Theory Phys.} \textbf{7} (2013), 145--223] as a special case, and provide supporting data for it. A~central role is played by the deformed Drinfel'd (or quantum) double of each umbral finite group~$G$, specified by a cohomology class in $H^3(G,U(1))$. We conjecture that in each of the 23 cases there exists a rule to assign an infinite-dimensional module for the deformed Drinfel'd double of the umbral finite group underlying the mock modular forms of umbral moonshine and generalised umbral moonshine. We also discuss the possible origin of the generalised umbral moonshine.}


\Classification{11F22; 11F37; 20C34}

\renewcommand{\thefootnote}{\arabic{footnote}}
\setcounter{footnote}{0}

\section{Introduction}\label{Intro}

Moonshine is a relation between finite groups and modular objects. The study of this relation started with the so-called monstrous moonshine \cite{conway_norton} and has recently been revived by the discovery of Mathieu moonshine \cite{Cheng2010_1,Eguchi2010a,Eguchi2010,Gaberdiel2010,Gaberdiel2010a} and umbral moonshine~\cite{UM,MUM}, which recovers Mathieu moonshine as a special case. While in moonshine one considers functions indexed by elements of a finite group, in generalised moonshine the relevant functions are indexed by pairs of commuting elements. The latter is a generalisation in that it recovers the former when the first element of the pair is set to be the identity element of the group. In this work we introduce {\em generalised umbral moonshine}, recovering the generalised Mathieu moonshine~\cite{Gaberdiel:2012gf} as a special case. With this we hope to contribute to the understanding of moonshine in the following two ways: 1)~by providing more examples of the moonshine relation, in particular working towards a complete analysis of umbral moonshine, and 2) by shedding light on the underlying structure of umbral moonshine which is still obscure.

The discovery of moonshine was initiated in the late 1970's with the realisation \cite{conway_norton} that there must be a representation $V^\natural_n$ of the monster group responsible for each of the Fourier coefficients of a set of distinguished modular functions. These modular functions~-- the Hauptmoduls $J_g$~-- are indexed by elements $g$ of the monster group:
\begin{gather*}
J_g =q^{-1}+\sum_{n=1}^\infty a_g(n)\,q^n, \qquad a_g(n) = {\rm Tr}_{V^\natural_n} g.
\end{gather*}
 Importantly, it was soon realised that this representation has an underlying algebraic structure: a vertex operator algebra (VOA), or a 2d chiral conformal field theory (CFT) in the physics language \cite{Bor_PNAS,FLM,FLMPNAS,FLMBerk}. In particular, the monster representation $V^\natural_n$ is nothing but the space of quantum states with energy (eigenvalue of ${\hat H}=L_0-c/24$) $n$ of a chiral CFT~-- the so-called monster CFT~-- that has monster symmetry.

Given a chiral conformal field theory with a discrete symmetry group $G$, one can also consider orbifolding the theory by subgroups of $G$. This leads to the concept of twisted sectors, which have symmetries given by elements of the group $G$ that commute with the twisting subgroup. Considering the action of this remaining symmetry group on the twisted sector quantum states leads to the {\it twisted-twined partition functions} indexed by two commuting elements~-- the twisting and the twining elements~-- of the original group $G$. At least heuristically, the existence of these twisted-twined partition functions is believed to be a conceptual partial explanation for Norton's {\it generalised monstrous moonshine} \cite{generalized_moonshine,generalized_moonshine2}, which we will review in more detail in Section~\ref{GMM}.

The landscape of moonshine has changed dramatically since the observation of Eguchi, Ooguri, and Tachikawa in 2010 \cite{Eguchi2010}. These authors pointed out a surprising empirical relation between the elliptic genus of $K3$ surfaces and the sporadic Mathieu group $M_{24}$. The ensuing development has led to the discovery of 23 cases of {\em umbral moonshine}, which recovers the above-mentioned $M_{24}$ relation as a special case. The important features of umbral moonshine include the following facts: 1) the relevant functions are so-called mock modular forms which have a~modified transformation property under the modular group, and 2) the 23 cases are organised in terms of 23 special (Niemeier) lattices of rank 24. In particular, their lattice symmetries dictate the finite groups featuring in umbral moonshine.

In order to understand umbral moonshine, it is crucial to know what the underlying algebraic structure is. Despite various recent advances~\cite{UMk3,UMLG,Duncan:2014tya,Harvey:2013mda,Harvey:2014cva,Kachru:2016ttg}, the nature of this structure is still obscure. An obvious first guess is that a 2d CFT will again be the relevant structure, especially given the CFT context in which the first case of umbral moonshine~-- the Mathieu moonshine~-- was uncovered. However, there is a salient incompatibility between the modular behaviour expected from a CFT partition function and that displayed in umbral moonshine that we will explain in Section~\ref{UMreview}.
As a result, new ingredients other than a conventional CFT/VOA are believed to be necessary for a uniform understanding of umbral moonshine. The quest for this new structure constitutes one of our motivations to study generalised umbral moonshine, which establishes the central role played by the deformed Drinfel'd double. Note that the importance of the third group cohomology in moonshine has been suggested by Terry Gannon and nicely demonstrated in~\cite{Gaberdiel:2012gf} in the context of generalised Mathieu moonshine.

The rest of the paper is organised as follows. To describe our results, we begin in Section~\ref{Background} by recalling the (generalised) monstrous moonshine, umbral moonshine, and reviewing the basic properties of the deformed Drinfel'd double of a finite group and its representations. In Section~\ref{GUM} we conjecture that there is a way to assign an infinite-dimensional module for the deformed Drinfel'd double of the umbral group underlying each of the 23 cases of generalised umbral moonshine and unpack this conjecture in a more explicit form in terms of the six conditions of generalised umbral moonshine. To provide evidence and explicit data for this conjecture, in Section~\ref{Umbral Groups, Cohomologies, and Representations} we present the group theoretic underpinnings of the generalised umbral moonshine and exposit our methodology and results. In Section~\ref{sec:mmf} we turn our attention to the modular side of moonshine and present our proposal for the twisted-twined functions of generalised umbral moonshine. In Section~\ref{Conclusions} we close the paper with a summary and discussions. To illustrate the content of the generalised umbral moonshine conjecture, in the appendix we provide the expansion of three examples of the $g$-twisted generalised moonshine functions belonging to three different lambencies, as well as the decompositions of the first few homogeneous components of the corresponding (conjectural) generalised umbral moonshine module into irreducible projective representations.

We also include two files accompanying this article; the first describes the explicit construction of the umbral groups that we employ, and the second describes the construction of the 3-cocycles.

\section{Background}\label{Background}
In this section we provide the relevant background for generalised umbral moonshine by reviewing the theories of generalised monstrous moonshine and umbral moonshine. We also review the basic properties of the deformed Drinfel'd double of a finite group and its representation theory.

\subsection{Generalised monstrous moonshine}\label{GMM}

The monster (or the Fischer--Griess monster) group $\mathbb M$ is the largest sporadic group. Monstrous moonshine connects the representation theory of the monster group and certain distinguished modular functions, namely the Hauptmoduls for genus zero subgroups of $\operatorname{PSL}_2(\RR)$.\footnote{A discrete subgroup $\Gamma\subset\operatorname{PSL}_2(\RR)$ is said to be genus zero if its fundamental domain on the upper-half plane, when suitabley compactified, is a genus zero Riemann surface. A function is said to be a Hauptmodul for $\G$ if it is an isomorphism from the compactified fundamental domain to the Riemann sphere.} In particular, the Fourier coefficients of the modular $J$-invariant are equated with the dimensions of the homogeneous components of a distinguished infinite-dimensional $\ZZ$-graded module for the monster group. This module was constructed (conjecturally) by Frenkel, Lepowsky, and Meurman \cite{FLMPNAS,FLMBerk}. It possesses the structure of a VOA and is later proven by Borcherds to be the monstrous module~\cite{borcherds_monstrous}. In this way, the $q$-powers of the $J$-invariant acquire the interpretation as the $L_0$-eigenvalues of the states in the corresponding chiral conformal field theory shifted by $-1$. Here $-1=-c/24$, where $c$ is the central charge of the Virasoro algebra.
More precisely, if we label by $V^\natural = \oplus_{n\geq -1}V^\natural_n$ the Hilbert space of quantum states of the monstrous chiral conformal field theory, we have
\begin{gather*}
 \Tr_{V^\natural} q^{L_0-c/24} : = \sum_{n\geq -1} q^n \dim V^\natural_n = J(\t) =q^{-1} + 196884 q+ 21493760 q^2 + \cdots,
\end{gather*}
where $q:= {\rm e}^{2\pi {\rm i} \t}$. More generally, for any $g\in \mathbb M$ the corresponding twined partition function, given by $ \Tr_{V^\natural} gq^{L_0-c/24} $, coincides with the Hauptmodul $J_g(\t)= q^{-1}+O(q)$ of a genus zero subgroup $\G_g$ of ${\rm PSL}_2(\RR)$.

An important way to construct new chiral CFTs is to ``orbifold" a chiral CFT that has a~discrete symmetry group~$G$~\cite{Dijkgraaf:1989hb,Dixon:1985jw,Dixon:1986jc, FLMBerk} by combining $g$-twisted modules for $g\in G$. The $g$-twisted sector Hilbert space, denoted by ${\cal H}_g$, has remaining symmetries given by elements of $G$ commuting with~$g$. By considering these symmetries, one can define the twisted-twined partition function
\begin{gather*}
Z_{(g,h)}(\t) = \Tr_{{\cal H}_g} h q^{L_0-c/24}
\end{gather*}
for all commuting pairs $(g,h)$ of $G$. Recall that in a (non-chiral) orbifold conformal field theory the above quantity is associated to a torus path integral interpretation with a $g$-twisted monodromy/defect for the spatial circle and $h$-boundary condition for the temporal circle. Hence, for any $\g =\sm abcd \in {\rm SL}_2(\ZZ)$ the transformation $\tau\mapsto {a\t+b\over c\t+d}$ of the torus complex structure is equivalent to the transformation $(g,h) \mapsto (g',h')= \big(g^a h^c, g^b h^d \big)$ of the twisting-twining elements. When working with a chiral theory, this consideration then leads to the expectation that $Z_{(g,h)}\big(\tfrac{a\t+b}{c\t+d}\big)$ and $ Z_{(g',h')}(\t) $ coincide up to a phase. Indeed, eight years after the monstrous moonshine conjecture, Norton proposed the following.

{\em Generalised monstrous moonshine conjecture} \cite{generalized_moonshine} (revised in~\cite{generalized_moonshine2}): There exists a rule that assigns to each element $g$ of the monster group $\mathbb M$ a graded projective representation $V(g)=\bigoplus_{n\in \QQ} V(g)_n$ of the centraliser group $C_{\mathbb M}(g)$, and to each pair $(g,h)$ of commuting elements of~$\mathbb M$ a~holomorphic function $T_{(g,h)}$ on the upper-half plane $\HH$, satisfying the following conditions:
\begin{enumerate}\itemsep=0pt
\item[\bf I] For any $\sm abcd\in {\rm SL}_2(\ZZ)$, $T_{(g^a h^c, g^b h^d)}(\t)$ is proportional to $T_{(g,h)}\big(\tfrac{a\t+b}{c\t+d}\big)$.
\item[\bf II] The function $T_{(g,h)}$ is either constant or a Hauptmodul for some genus zero congruence group.
\item[\bf III] The function $T_{(g,h)}$ is invariant up to constant multiplication under simultaneous conjugation of the pair $(g,h)$ in ${\mathbb M}$.
\item[\bf IV] There is some lift $\til h$ of $h$ to a linear transformation on $V(g)$ such that
\begin{gather*} T_{(g,h)}(\t) = \sum_{n\in \QQ}\operatorname{Tr}_{V(g)_n} \til h q^{n-1}.
\end{gather*}
\item[\bf V] The function $T_{(g,h)}$ coincides with the $J$-invariant if and only if $g=h=e\in \mathbb M$.
\end{enumerate}

The proof of this conjecture has been recently announced \cite{carnahan12}, building on previous work \cite{Carnahan2008,MR2904095,carnahan_miyamoto,Dong2000,Hoehn03}. See also~\cite{MR1372731} for a previous result and \cite{carnahan12} for a more in-depth review of the literature. With the crucial exception of the genus zero property ({\bf II}), all the above properties of generalised moonshine can be understood in the physical framework of holomorphic orbifolds along the lines discussed above \cite{Dixon:1988qd}. See~\cite{Paquette:2016xoo} for a recently proposed physical interpretation of the genus zero property.

\subsection{Umbral moonshine}\label{UMreview}

Historically, the study of umbral moonshine started with the discovery of Mathieu moonshine, initiated by a remarkable observation made in the context of $K3$ conformal field theories in~\cite{Eguchi2010}. The Mathieu moonshine can be phrased as the fact that the Fourier coefficients of a set of mock modular forms can be equated with the group characters of a distinguished infinite-dimensional $\ZZ$-graded module for $M_{24}$. Later it was realised in \cite{UM,MUM} that the Mathieu moonshine is but an instance of a larger structure, which was named umbral moonshine. Umbral moonshine is labelled by Niemeier lattices, the unique (up to isomorphism) 23 even self-dual positive-definite lattices in 24 dimensions with a non-trivial root system. Recall that the root system is the sub-lattice generated by the root vectors, the lattice vectors of norm squared 2. These 23 Niemeier lattices are uniquely specified by their root systems, which are precisely the 23 unions of simply-laced (ADE) root systems with the same Coxeter number whose ranks are 24. The relevant finite groups, the so-called umbral groups, are given by the automorphism of the corresponding lattice quotiented by the Weyl reflections with respect to the lattice root vectors.

It turns out to be natural to associate a~genus zero subgroup of ${\rm SL}_2(\ZZ)$ of the form $\Gamma_0(m)+e, f, \dots $ to each of the 23 Niemeier lattices \cite{omjt, MUM}, where $\Gamma_0(m)+e, f, \dots $ denotes the subgroup of ${\rm SL}_2(\RR)$ obtained by attaching the corresponding Atkin--Lenher involutions $W_e, W_f, \dots $ to the congruence subgroup~$\G_0(m)$. Following \cite{MUM} we use the shorthand $\ell = m +e, f, \dots$ for the genus zero group $\Gamma_0(m)+e, f, \dots $ and call $\ell$ the lambency of the corresponding instance of umbral moonshine. The level $m$ of the group $\Gamma_0(m)+e, f, \dots $ coincides with the index of the corresponding mock Jacobi forms, so we say that $m$ is the index of that lambency. We also observe that $m$ is the Coxeter number of the root system of the corresponding Niemeier lattice:
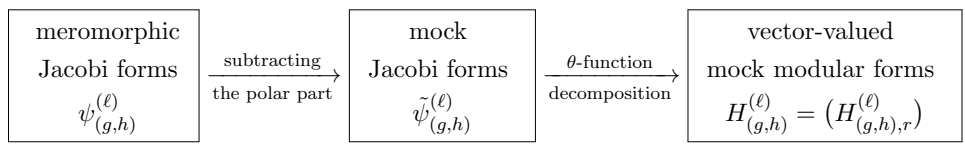
\begin{gather*}\label{eq:modular_objects}\resizebox{0.8\hsize}{!}{%
 $\boxed{ \begin{array}{c}
\text{ meromorphic}\\
\text{ Jacobi forms}\\~\psi^{(\ell)}_{(g,h)}\end{array}} \xrightarrow[\text{the polar part}]{\text{subtracting}}
\boxed{\begin{array}{c}
\text{mock }\\
\text{Jacobi forms }\\\tilde\psi^{(\ell)}_{(g,h)}\end{array}} \xrightarrow[\text{decomposition}]{\theta\text{-function}}
\boxed{\begin{array}{c}
\text{vector-valued }\\ \text{mock modular forms }\\
H^{(\ell)}_{(g,h)}=\big(H^{(\ell)}_{(g,h),r}\big) \end{array} }$
}\end{gather*}

There are a few ways to view the modular objects of umbral moonshine, which are summarised in diagram~(\ref{eq:modular_objects}) and which we will now explain. We start by defining a number of actions on functions defined on ${\mathbb H}\times{\mathbb C}$ that will be used in the modularity criteria. For a given $m\in \ZZ_+$ we defined the index $m$ (generalised) elliptic action by\footnote{Note that usually the elliptic action is only defined for $(\l,\m)\in \ZZ^2$, in which case the factor $\ex(m\l\m)$ is unity and is not included in the standard definition.}
\begin{gather}\label{def:ellip}
(\f\big|_m(\lambda,\mu))(\tau,\zeta):=\f(\t, \z+\l \t +\m)\ex\big( m\big(\l^2 \t + 2\l \z + \l\m\big)\big)
\end{gather} for any $(\l,\m)\in \RR^2$. Here and in the rest of the paper we will use the shorthand notation $\ex(x):={\rm e}^{2\pi {\rm i} x}$. Moreover, for $2k\in \ZZ$ and $\gamma =\sm abcd \in {\rm SL}_2(\ZZ)$ we define the modular transformation
\begin{align}
\label{modular}
(\f\big\lvert_{k,m}\g )(\tau,\zeta) &:=
\f\left(\frac{a\t+b}{c\t+d},\frac{\z}{c\t+d}\right)
\frac1{(c\t+d)^{{k}}}
\ex\left(- \frac{c m\z^2}{c\t+d}\right).
\end{align}
It is easy to check that the elliptic and modular operators satisfy
\be\label{modularelliptic}
\f \lvert_m(x,y) \lvert_{k,m} \g = \f \lvert_{k,m} \g \lvert_m(\g(x,y)).
\ee
for all $\g\in{\rm SL}_2(\ZZ)$ and $(x,y)\in\RR^2$, where we have defined $\gamma(x,y) := (ax+cy, bx+dy)$.

Associated to each lambency $\ell=m +e, f, \dots$ and each element $g$ of the corresponding umbral group $G^{(\ell)}$ is a so-called meromorphic Jacobi form $\psi^{(\ell)}_{g}$. In particular, $\psi^{(\ell)}_{g}$ transforms as a weight one index $m$ Jacobi form for a congruence subgroup $\Gamma_g\subseteq {\rm SL}_2(\ZZ)$, possibly with a non-trivial multiplier (cf.\ Sections~\ref{GUM} and~\ref{New Groups}). Namely, we have
\begin{gather*}
\psi^{(\ell)}_{g}\big\lvert_m(\lambda,\mu) =\psi^{(\ell)}_{g}\end{gather*}
 and
 \begin{gather}\label{untwisted_modu}
\varsigma_g(\gamma) \psi^{(\ell)}_{g}\big\lvert_{1,m} \g =\psi^{(\ell)}_{g}
 \end{gather}
for all $(\l,\m)\in\ZZ^2$, $\g\in \Gamma_g$, and for some $\varsigma_g\colon \Gamma_g\to \CC^\times$. Moreover, $\psi^{(\ell)}_{g}$ is meromorphic as a~function $\zeta \mapsto \psi_{g}^{(\ell)}(\t,\zeta)$ and can have poles at most at $2m$-torsion points, ${\zeta\in{1\over 2m}\ZZ+{\t\over 2m}\ZZ}$.

As explained in \cite{Dabholkar:2012nd, Zwegers2008}, meromorphic Jacobi forms such as $\psi_{g}^{(\ell)}$ lead to vector-valued mock modular forms. Recall that a holomorphic function $f\colon \HH\to \CC$ is said to be a~weight $k$ mock modular form for $\G\subseteq {\rm SL}_2(\ZZ)$ with shadow $g$ if $g$ is a holomorphic modular form of weight $2-k$ such that the non-holomorphic function
\begin{gather*}
\hat f(\t):=f(\t)+(4{\rm i})^{k-1}\int_{-\bar \t}^\infty (z+\tau)^{-k} \overline{g(-\bar \t)} {\rm d}z
\end{gather*}
transforms as a weight $k$ modular form for $\G$. To explain the relation between meromorphic Jacobi forms and vector-valued mock modular forms in the current context, we first define the holomorphic function $\tilde \psi_{g}^{(\ell)}$ from the meromorphic function $\psi_{g}^{(\ell)}$ by subtracting its ``polar part''
\begin{gather}\label{polarpart}
\tilde\psi_{g}^{(\ell)}(\tau,\zeta) := \psi_{g}^{(\ell)}(\tau,\zeta) - \psi_{g}^{{(\ell)},P}(\tau,\zeta) .
\end{gather}
 The polar part $\psi_{g}^{{(\ell)},P }$ is
 \begin{gather}\label{polar_A_type}
 \psi_{g}^{{(\ell)},P } =\sum_{a,b\in \ZZ/2m\ZZ}\chi^{(a,b)}_{g} \mu_{m}\lvert_m\big(\tfrac{a}{2m},\tfrac{b}{2m}\big), \end{gather}
 with $\chi^{(a,b)}_{g}\in \CC$ given by the characters of certain representations of the umbral group $\Gell$~\cite{MUM}. In the above we have used the generalised Appell--Lerch sum
 \begin{gather*}\label{gAPsum}
\m_{m} (\tau,\zeta) = - \sum_{k\in \ZZ} q^{m k^2} y^{2 m k} \frac{1+ yq^k }{ 1-yq^k},\end{gather*}
with $q=\ex(\t)$ and $y=\ex(\z)$.

It follows from the mock modular property of the Appell--Lerch sum (see, for instance, \cite{Dabholkar:2012nd, Zwegers2008}) that $\tilde \psi_{g}^{(\ell)}$ is an example of the so-called weak mock Jacobi forms.
Moreover, from the elliptic transformations of $ \psi_{g}^{(\ell)}$ (cf.\ equation~\eqref{elliptic_trans}) and $\m_{m}$ we see that $\tilde \psi_{g}^{(\ell)}$ satisfies the same elliptic transformation as $ \psi_{g}^{(\ell)}$, and therefore admits the following theta function decomposition~\cite{Dabholkar:2012nd}
\begin{gather}\label{mockjacobi}
\tilde \psi_{g}^{(\ell)}(\tau,\zeta) = \sum_{r\in \ZZ/2m\ZZ} H_{g,r}(\t)\theta_{m,r}(\tau,\zeta) ,
\end{gather}
where
\begin{gather}\label{def:theta1}\th_{m,r}(\tau,\zeta) = \sum_{k=r\xmod{2m}} q^{k^2/4m} y^k.
\end{gather}

Recall that, when regarded as a vector-valued function with $2m$ components, $\th_m=(\th_{m,r})$ transforms under $\gamma\in {\rm SL}_2(\ZZ)$ as
\begin{gather}\label{theta_mult}
\rho_{m}(\gamma) \th_m\lvert_{{1\over 2},m}\gamma= \th_m
\end{gather}
with the multiplier system $ \rho_{m}\colon {\rm SL}_2(\ZZ) \to {\rm GL}(2m,\CC)$ whose detailed description can be found in, for instance, \cite[Section~5]{eichler_zagier}. From the modular property of $\psi_{g}^{(\ell)}$ \eq{untwisted_modu} and the mock modular property of the Appell--Lerch sum \cite{Dabholkar:2012nd,Zwegers2008}), we see that $H_{g} = ( H_{g,r})$ is a weight $1/2$ vector-valued mock modular form\footnote{In general $H_{g,r}(\t)$ grows exponentially near at most one cusp of the group $\Gamma_{g}$. Hence, it is {\em weakly holomorphic} vector-valued mock modular forms that we encounter in generalised umbral moonshine. To avoid clutter we will skip the adjective ``weakly holomorphic'' in the remaining of the paper.} for a certain $\Gamma_{g}\subseteq {\rm SL}_2(\ZZ)$ that we will define in \eq{def:gamma_gh}.

In the rest of the paper, we will generalise the above relation between meromorphic Jacobi forms, mock Jacobi forms, and vector-valued mock modular forms, to functions labelled by commuting pairs of elements of $G^{(\ell)}$. This relation between $\psi_{(g,h)}$, $\til \psi_{(g,h)}$ and $ H_{(g,h)}$ is summarised in diagram~(\ref{eq:modular_objects}).

Just like the Hauptmoduls in monstrous moonshine, the mock modular forms $H^{(\ell)}_g$ are also special as we shall explain now. Recall that Cappelli, Itzykson and Zuber found an ADE classification of modular invariant combinations of~$A_1$ affine characters~\cite{Cappelli:1987xt}. This classification, combined with the ADE data from the root systems of the Niemeier lattices, leads to a specification of the modularity property of the mock modular forms of umbral moonshine as described in~\cite{MUM}. The umbral moonshine conjecture states that these vector-valued mock modular forms~$H^{(\ell)}_{g}$ for $g=e$, are in fact the unique such forms with the slowest possible growth in their Fourier coefficients given their modularity properties. In particular, $H_e^{(\ell)}$ coincides with the Rademacher sum specified by its modular property and its pole at the cusps~\cite{MUM}. See \cite{Cheng2011,2012arXiv1212.0906C,mumcorr,Rad_whalen} for a~discussion on the cases where~$g\neq e$.

After specifying the finite group $G^{(\ell)}$ and the vector-valued mock modular forms $H^{(\ell)}_{g}=(H_{g,r}^{(\ell)})$ for all elements $g \in G^{(\ell)}$, the umbral moonshine conjecture states that there exists a $G^{(\ell)}$-module whose graded characters coincide with the Fourier expansion of $H^{(\ell)}_{g}$. Namely, there exists an infinite-dimensional $\Gell$-module
\begin{gather*}
K^{(\ell)} =\bigoplus_{\substack{r\in \{1,\dots,m-1\}\\ \a\in \QQ_{>0}}} K^{(\ell)}_{r,\a}
\end{gather*}
 such that for all $g\in G$ and $r\in \{1,\dots,m-1\}$ we have
\begin{gather}\label{um_conj}
 H^{(\ell)}_{g,r}(\t)= -2 \delta_{r^2,1}^{[4m]} q^{-1/4m} + \sum_{\a\in \QQ_{>0}} q^\a \text{Tr}_{K^{(\ell)}_{r,\a}}g.
\end{gather}
In the above equation we defined
\begin{gather*}
\delta_{a,b}^{[C]} = \begin{cases} 1 &{\rm if}~a=b\xmod{C}, \\ 0 &{\rm otherwise}.\end{cases}
\end{gather*}
Corresponding to the Niemeier lattice with root system $A_1^{24}$ is the case of umbral moonshine with $\ell =2$ and $\Gell\cong M_{24}$. In this case, the vector-valued mock modular forms have just a single independent component since $H^{(2)}_{g,r}=\pm H^{(2)}_{g,1}$ for $r=\pm 1\xmod{4}$ and $H^{(2)}_{g,r}=0$ for $r=0\xmod{2}$, and the above conjecture recovers that of Mathieu moonshine.

The above conjecture was proven in \cite{proof, Gannon:2012ck}, in the sense that the existence of the module~$K^{(\ell)}$ has been established using properties of (mock) modular forms. However, among the 23 cases of umbral moonshine, modules have only been constructed for the eight simpler cases \cite{Anagiannis:2017src,Cheng:2017grj, Duncan:2014tya,Duncan:2017bhh}. A uniform construction of the umbral moonshine modules is to the best of our knowledge not yet in sight and is expected to be the key to a true understanding of this new moonshine phenomenon. One of the most puzzling features in this regard is the non-vanishing weight of the umbral moonshine function: if one assumes that the umbral moonshine module, just like the monstrous moonshine module, possesses the structure of a conventional chiral CFT, than the usual physics lore mentioned in the last subsection as well as mathematical results on VOAs~\cite{Zhu_ModInv} lead to the expectation that one should be able to attach a weight zero modular object to the moonshine module. This process can be done in a straightforward manner for the eight cases of umbral moonshine which correspond to the eight Niemeier lattices with $A$-type root systems, but not for the remaining fifteen cases involving $D$- and/or $E$-type root systems, at least not if we require the modular object to be holomorphic as usual. Moreover, among all the A-type cases, it was argued in~\cite{UM} that apart from the~$A_1^{24}$ case of Mathieu moonshine, the associated weight zero Jacobi forms cannot coincide with a partition function or elliptic genus of a ``usual'' (chiral) conformal field theory. This is due to the expectation that the NS-NS ground states lead to a non-vanishing $q^0 y^{m-1}$ term in the index $m$ weight zero Jacobi form, which is absent in all cases other than the~$A_1^{24}$ case. In this regard, the $A_1^{24}$ case of umbral moonshine is set apart from the other cases that the mock modular forms can be regarded as arising from the supersymmetric index of a physical conformal field theory~-- any $K3$ non-linear sigma model in fact. This more involved modular property suggests that a novel approach is likely to be needed in order to understand umbral moonshine in general.

\subsection{Deformed Drinfel'd double and its representations}\label{Deformed Quantum Double and its Representations}

To set the stage for the generalised umbral moonshine conjecture, we will quickly review the definition of a deformed Drinfel'd (or quantum) double of a finite group and its representations. Drinfel'd introduced the quantum double construction, which associates to a~Hopf algebra~$\mathcal A$ a~quasi-triangular Hopf algebra~$D({\cal A})$. This quasi-triangular Hopf algebra contains ${\cal A}$ as well as its dual algebra ${\cal A}^\ast$ as Hopf algebras and $D({\cal A})= {\cal A}^\ast \otimes A $ as a vector space~\cite{MR934283}. Let~$G$ be a~finite group and ${\cal A}=kG$ be the group algebra over field $k$. It was argued that the corresponding quantum double for the case $k=\CC$, often denoted simply by $D(G)$, is relevant for 3d TQFTs~\cite{Dijkgraaf:1989pz} and relatedly 2d orbifold CFTs \cite{Roche:1990hs, Dijkgraaf:1989hb}. The finite group~$G$ plays the role of the discrete gauge group in te former case and that of the orbifold group in the latter case. Here and in the rest of the paper we will work with $k=\CC$ exclusively.

In \cite{Roche:1990hs}, Dijkgraaf, Pasquier and Roche introduced an interesting generalisation of this construction. For $G$ a finite group and $M$ a $G$-module, one can define the group cohomology $H^n(G,M):= \operatorname{Ker} \pa^{(n)} / \operatorname{Im} \pa^{(n-1)}$ where the coboundary operator $\pa^{(n)}$ on the space of $n$-cochains is given by
\begin{gather*}
\big(\pa^{(n)} \l\big)(h_1,\dots,x_n,x_{n+1}) = \l(x_2,\dots,x_{n+1}) + {(-1)^{n+1}} \l(h_1,\dots,x_{n})\\
\hphantom{\big(\pa^{(n)} \l\big)(h_1,\dots,x_n,x_{n+1}) =}{} +\sum_{i=1}^n {(-1)^i} \l (h_1,\dots,x_i x_{i+1} ,\dots, x_{n+1}).
\end{gather*}
From now on we take $M=U(1)$ (on which $G$ acts trivially), and to avoid an overload of notation we will often conflate $x\in \RR/\ZZ\cong U(1)$ with ${\rm e}^{2\pi {\rm i} x} \in \CC^\times$ when the meaning is clear from the context. Given $G$ and a 3-cocycle $\omega\colon G\times G\times G \to \CC^\times$, the deformed quantum double~$D^\omega(G)$ is a quasi-triangular quasi-Hopf algebra generated by the elements\footnote{Here we use the notation in \cite{Bantay:1990tr}. Some other common choices of notation for $P(x)Q(y)$ include $\d_xy$, $\mbox{\fontsize{6}{10}\selectfont $x$}\,\underset{y}{\halfbox{ \textcolor{white}{X}}}$ and $\langle \xleftarrow[]{y} x \rangle$. We will work only with the so-called normalised cochains, namely functions $u\colon G^n \to \CC^\times$ satisfying $u(x_1,\dots,x_n)=1$ whenever at least one of the $n$ group elements $x_1,\dots,x_n$ coincides with the identity element~$e$. This corresponds to requiring $\rho(e)={\mathbf 1}_V \in \operatorname{End}(V)$ in the $G$-representation $(\rho,V)$.}~$P(x)$ and~$Q(y)$, with \mbox{$x,y \in G$}~\cite{Roche:1990hs}. The multiplications are given by
\begin{gather*}
P(x)P(y) = \begin{cases} P(x) &{\rm if }~x=y, \\ 0 &\text{otherwise},\end{cases}\qquad
Q(x) Q(y) = \th(x,y) Q(xy),
\end{gather*}
while the co-multiplication are given by
\begin{gather*}
\D P(x) = \sum_{z\in G} P(z) \otimes P\big(z^{-1}x\big), \qquad \D Q(x) = \eta(x) Q(x) \otimes Q(x).
\end{gather*}
The multiplication and co-multiplication rules are determined by the {3-cocycle} $\omega$ as
\begin{gather*}
\th(x,y)= \sum_{z\in G} \th_z(x,y) P(z), \qquad
\eta(z) = \sum_{x,y\in G} \eta_z(x,y) P(x)\otimes P(y),
\end{gather*}
where
\begin{gather} \label{def:2cycle}
\theta_g(x,y) := \frac{\omega(x,g\lvert^{x},y)}{\omega(g,x,y)\omega(x,y,g\lvert^{xy})}, \\ \label{def:2cycle2}
\eta_g(x,y) := \frac{\omega(x,g,y\lvert^{x})}{\omega(x,y,g)\omega(g,x\lvert^g,y\lvert^g)}.
\end{gather}

Moreover, we have
\begin{gather*}
Q(x)^{-1} P(y) Q(x) = P\big(y\lvert^x\big),
\end{gather*}
where we use the notation $y\lvert^x := x^{-1}y x$ to denote the conjugate of $y$ by $x$. The identity element of $D^\omega(G)$ is $Q(1) = \sum\limits_{z\in G} P(z)$.

Note that $\theta_g(x,y)=\eta_g(x,y)$ when $x, y \in C_G(g)$, the centralizer of $g$ in $G$. It is easy to check that the restriction $c_g:=\theta_g\lvert_{C_G(g)\times C_G(g)}=\eta_g\lvert_{C_G(g)\times C_G(g)} $ is a 2-cocycle and moreover, that the corresponding map $Z^3(G,U(1)) \to Z^2(C_G(g),U(1)) $ induces a map $H^3(G,U(1))\to H^2(C_G(g),U(1))$.

The fact that the above defines a quasi-Hopf algebra is guaranteed by the 3-cocycle condition $\pa^{(3)} \omega=0$. In particular, the co-multiplication, instead of satisfying $({\mathbf 1}\otimes \D)\cdot \D = (\D\otimes {\mathbf 1})\cdot \D$ as in the case of a Hopf algebra, satisfies
\begin{gather*}
({\mathbf 1}\otimes \D)\cdot \D (a) =\varphi (\D\otimes {\mathbf 1})\cdot \D(a) \varphi^{-1} \qquad \forall a \in D^{\o}(G)
\end{gather*}
in a quasi-Hopf algebra, where the intertwiner is given by
\begin{gather*}
\varphi = \sum_{x,y,z\in G} \omega(x,y,z) P(x)\otimes P(y)\otimes P(z).
\end{gather*}
It is easy to see that changing $\w$ without changing the cohomology class leads to an isomorphic quasi-Hopf algebra \cite{Roche:1990hs}. In the language of holomorphic orbifold the non-triviality of the 3-cocycle signals the failure of the fusion between operators belonging to different twisted sectors to be associative. See~\cite{Roche:1990hs} and references therein for more details on $D^{\o}(G)$ including the $R$-matrices and the antipode.

Next we will summarize the basic properties of the representations of $D^\o(G)$. See, for instance, \cite{Altschuler:2003nn,Bantay:1990tr, Roche:1990hs} for more information. Firstly, we will establish that given an $\o$-compatible projective representation of $C_G(g)$ for some $g\in G$ we can build a corresponding $D^\o(G)$-representation via the so-called DPR induction.

Recall that a projective representation $(\rho,V)$, where $\rho\colon \Gg\to \operatorname{End}(V)$, satisfies
 \begin{gather*}
 \rho(h_1)\rho(h_2) = c^\rho(h_1,h_2)\,\rho(h_1h_2)
 \end{gather*}
for any $h_1,h_2\in \Gg$, where $c^\rho\in Z^2(\Gg,U(1))$. Given a ${g_A}\in G$ and a projective representation $(\rho_{g_A}, V_{g_A})$ of $C_G(g_A)$ corresponding to the 2-cocycle $\th_{g_A}\in Z^2(C_G(g_A),U(1))$, we can construct a {\em DPR-induced} representation of $D^\o(G)$ in the following way~\cite{Roche:1990hs}. Let ${\cal B}_{g_A}$ be the subalgebra of $D^\o(G)$ spanned by elements of the form $P(g)Q(x)$ with $g\in G$ and $x\in C_G(g_A)$. We define the action of ${\cal B}_{g_A}$ on $V_{g_A}$ by $\p(P(g)Q(x))= \d_{g,g_A} \rho_{g_A}(x)$. The DPR-induced $D^\o(G)$ representation, denoted $\big(\pi_{(\rho_{g_A}, V_{g_A})}, {\rm Ind}^{{\rm DPR}}(V_{g_A})\big)$, is given by ${ \rm Ind}^{{\rm{DPR}}}(V_{g_A}):= \CC[G]\otimes_{{\cal B}_{g_A}} V_{g_A}$, where we identify~$\CC[G]$ as the subalgebra spanned by $\sum\limits_{z\in G}P(z) Q(x)$ with $x\in G$.
 Explicitly, to describe the action of $D^\o(G)$ we choose a set of representatives $\{x_1,\dots,x_n\}$ of $G/C_G(g_A)$, and it can be checked that, for $v\in V_{g_A}$,
 \begin{gather*}
\pi_{(\rho_{g_A}, V_{g_A})}(P(g)Q(x))(x_j\otimes v) ={\th_{g_A}(x,x_j)\over\th_{g_A}(x_k,h)} \d_{g\lvert^x, g_j} (x_k\otimes \rho_{g_A}(h) v),
\end{gather*}
where $h\in C_G(g_A)$ is determined by $xx_j=x_k h$.

Secondly, any irreducible representation of $D^\o(G)$ is labelled by a conjugacy class $A$ of $G$ and an irreducible projective representation of $C_G(g_A)$ with the 2-cocycle $\th_{g_A}$, where $g_A\in A$. Moreover, such a~irreducible representation of $D^\o(G)$ is equivalent to the corresponding DPR-induced representation described above. As a result, studying the representations of $D^\o(G)$ is equivalent to studying the projective representations of all the centraliser subgroups with the 2-cocycles given by $\o$ via~\eq{def:2cycle}, a fact that we will exploit in order to make explicit conjectures in the next subsection.

The representations of a deformed quantum double share many features with those of finite groups. For instance there is an orthogonality relation among irreducible representations. For a given representation $V$ of $D^\o(G)$, we write its character as $\operatorname{Ch}_{(x,y)}^V : =\Tr_V (P(x)Q(y))$. These characters satisfy the following properties,
\begin{gather}
\operatorname{Ch}_{(x,y)}^V = 0 \qquad {\rm unless}\quad xy = yx,\nonumber \\ \label{conjugation}
\operatorname{Ch}_{{(x\lvert^z,y\lvert^z)}}^V = \frac{\th_x(z,y\lvert^z)}{\th_x(y,z)}\operatorname{Ch}_{(x,y)}^V ,
\end{gather}
which can also be checked explicitly for the DPR-induced representations. Moreover, $\rm{Rep}(D^\o\!(G))\!$ is a modular category. In terms of the characters, the modular group acts as
\begin{gather}\label{Tmod}
T\operatorname{Ch}_{(x,y)}^V = \th_x(x,y) \operatorname{Ch}_{(x,xy)}^V, \\ \label{Smod}
S\operatorname{Ch}_{(x,y)}^V ={1\over \th_{y}(x,x^{-1})} \operatorname{Ch}_{(y,x^{-1})}^V ,
\end{gather}
satisfying $S^2 = (ST)^3$, $S^4={\rm id}$. In particular, the fusion rules between different representations satisfy the Verlinde formula given by the above $S$-matrix.

\section{Generalised umbral moonshine}\label{GUM}

The {\em generalised umbral moonshine conjecture} states that for each of the 23 lambencies of umbral moonshine with umbral group $G$ there is a 3-cocycle $\o\colon G\times G \times G\to \CC^\times$ such that there exists
an infinite-dimensional module for $D^\o(G)$ underlying a set of distinguished meromorphic Jacobi forms that we will specify in Section~\ref{sec:mmf}. In the following we will re-formulate the conjecture in explicit terms, using the relation to the projective representations of the centraliser subgroups discussed in Section~\ref{Deformed Quantum Double and its Representations}. This leads to the first five of the six conditions that we will discuss in this section. The sixth condition is the counterpart of the Hauptmodul condition and demonstrates the special property of the generalised umbral moonshine functions.

Explicitly, given the lambency $\ell$ with index $m$ and the corresponding umbral group $G$, we propose in Section~\ref{sec:mmf} the twisted-twined functions $\psi_{(g,h)}$ (and the corresponding twisted-twined mock modular forms $H_{(g,h)}$) for each commuting pair $(g,h)$. We conjecture that they satisfy the six conditions listed below. The relation between $\psi_{(g,h)}$ and the mock modular forms $H_{(g,h)}$ is as described in~(\ref{eq:modular_objects}).
\begin{enumerate}\itemsep=0pt
\item[\bf I]{\bf Modularity.} The modularity condition first states that the function $\psi_{(g,h)}$, which is meromorphic as a~function $\zeta \mapsto \psi_{(g,h)}(\t,\zeta)$ and can have poles (at most) at $2m$-torsion points, $\zeta\in{1\over 2m}\ZZ+{\t\over 2m}\ZZ$, transforms in the same way as an index $m$ Jacobi form under the elliptic transformation (cf.\ equation~\eqref{def:ellip}). Namely, for all $(\l,\m)\in\ZZ^2$,
\begin{empheq}[box=\widefbox]{equation}\label{elliptic_trans}
\psi_{(g,h)}\big\lvert_m(\lambda,\mu) =\psi_{(g,h)}.\end{empheq}

Second, under the weight one modular transformation they satisfy for any ${g,h\in G}$ and any $\g\in {\rm SL}_2(\ZZ)$
\begin{empheq}[box=\widefbox]{equation}\label{sl2z}
\varsigma_{(g,h)}(\g) \psi_{(g,h)}\big\lvert_{1,m}\g = \psi_{\gamma(g,h)}
\end{empheq}
 for some ${\varsigma_{(g,h)}\colon {\rm SL}_2(\ZZ) \to \CC^\times}$. In the above we have defined $\gamma(g,h):= \big(g^a h^c , g^b h^d\big)$ for two commuting elements $g,h\in G$.

This condition reflects the modular properties (\ref{Tmod}), (\ref{Smod}) of $D^\o(G)$ representations, embodied by the meromorphic Jacobi forms $\psi_{(g,h)}$, and can be regarded as the counterpart of condition {\bf I} of generalised monstrous moonshine.

\item[\bf II]{\bf 3-cocycle.} This condition states that there exists a 3-cocycle $\omega\colon G\times G\times G \to \CC^\times$ compatible with the multiplier system $\varsigma_{(g,h)}$ in \eq{sl2z} in the following way.
For any $g\in G$ let $\th_g$ be defined as in~\eq{def:2cycle}. Then the compatibility conditions read
\begin{subequations}\label{def:multiplier}
\begin{empheq}[box=\widefbox]{gather}
\varsigma_{(g,h)}\left(\sm1{1}01\right) = \theta_g(g,h)^{-1},\\
\varsigma_{(g,h)}\left(\sm0{-1}10\right) = \theta_h\big(g,g^{-1}\big),\end{empheq}
\end{subequations}
which determine $\varsigma_{(g,h)}$ in terms of $\o$ together with \eq{sl2z}.

This condition captures the modular properties (\ref{Tmod}), (\ref{Smod}) of representations of the quantum deformed double~$D^\o(G)$. It is believed that a 3-cocycle underlies the generalised monstrous moonshine in a similar way \cite{Gaberdiel:2012gf}, although it was not mentioned in any of the five stated conditions of generalised monstrous moonshine.

\item[\bf III]{\bf Projective class function.} This condition states that $\psi_{(g,h)}$ is a projective class function, satisfying for all $k\in G$
\begin{empheq}[box=\widefbox]{equation}\label{proj_class_function}
\psi_{(g\lvert^k,h\lvert^k)}=\xi_{(g,h)}(k) \psi_{(g,h)}
\end{empheq}
 with the phases $\xi_{(g,h)}(k)$ determined by the 3-cocycle $\o$ as
 \begin{gather*}
\xi_{(g,h)}(k)= \frac{\th_g(k,h\lvert^k)}{\th_g(h,k)}.
 \end{gather*}

This property together with the modularity property \eq{sl2z} implies that the ${\rm SL}_2(\ZZ)$ subgroup for which $H_{(g,h)}$ is a vector-valued mock modular form is given by
\begin{gather}\label{def:gamma_gh}
\Gamma_{g,h} := \big\{ \gamma\in {\rm SL}_2(\ZZ) \big\lvert \gamma(g,h) = \big(g\lvert^k,h\lvert^k\big) \text{ for some } k\in G\big\}.
\end{gather}

This condition captures the property~(\ref{conjugation}) of representations of the quantum deformed double~$D^\o(G)$ and is the counterpart of condition~{\bf III} of generalised monstrous moonshine.

\item[\bf IV]{\bf Finite group.} The vector-valued mock modular forms $H_{(g,h)}=(H_{(g,h),r})$ have the following moonshine relation to the finite group $G=G^{(\ell)}$. Namely, there is a way to assign an infinite-dimensional projective $C_G(g)$-module
\begin{gather*}
K^{g} =\bigoplus_{\substack{r\in \ZZ/2m\ZZ\\ \a\in \QQ_{>0}}} K^{g}_{r,\a}
\end{gather*}
with the 2-cocycle given by $c^\rho=\th_g$, such that for all $g\in G$, $h\in C_G(g)$ and $r\in \ZZ/2m\ZZ$,
\begin{empheq}[box=\widefbox]{equation}\label{finite_grp_property}
 H_{(g,h),r}(\t)= c_{(g,h),r}q^{-1/4m} + \sum_{\a\in \QQ_{>0}} q^\a \text{Tr}_{K^{g}_{r,\a}}h .
\end{empheq}
 The coefficient $c_{(g,h),r}$ equals $\mp 2$ whenever $H_{(g,h),r}(\t) = \pm H_{(e,h),1}(\t)$ and vanishes otherwise.

The above condition describes the $D^\o(G)$-module as DPR-induced modules as explained in Section~\ref{Deformed Quantum Double and its Representations} and can be regarded as the counterpart of condition {\bf IV} of generalised monstrous moonshine.

\item[\bf V]{\bf Consistency.} This condition states that the generalised umbral moonshine is compatible with umbral moonshine in the sense that $\psi_{(e,h)}^{(\ell)}$ coincides with the weight one meromorphic Jacobi form~$\psi^{(\ell)}_{h}$ given in~\cite{UM,MUM}. In particular, the resulting mock modular form $H_{(e,h)}$ coincides with the McKay--Thompson series $H^{(\ell)}_{h}=(H^{(\ell)}_{h,r})$ proposed in~\cite{UM,MUM}.

This condition establishes the relation between generalised umbral moonshine and the original umbral moonshine, and can be regarded as the counterpart of condition~{\bf V} of generalised monstrous moonshine.

\item[\bf VI]{\bf Pole structure.} This condition states that the mock modular form $H_{(g,h)}$ is bounded at all but at most one orbit of the cusps $\QQ\cup \{{\rm i}\inf\}$ under $\G_{(g,h)}$. It serves as the counterpart of condition {\bf II} (the Hauptmodul condition) of generalised monstrous moonshine.
\end{enumerate}

An immediate consequence of the generalised umbral moonshine conjecture is certain non-trivial constraints on their twisted-twined functions. To illustrate this, let us consider a pair $(g,h)$ of commuting elements in~$G$ that is conjugate to $\big(g^{-1},h^{-1}\big)$. The conditions~\eq{sl2z} and~\eq{proj_class_function} lead to the equality
\begin{gather}\label{eq:obstructions}
\varsigma_{(g,h)}\left(\sm{-1}00{-1}\right) \psi_{(g,h)}= \xi_{(g,h)}(k)\psi_{(g,h)},
\end{gather}
where $k$ a group element satisfying $\big(g^{-1},h^{-1}\big)=\big({g}\lvert^k,h\lvert^k\big)$. It follows that $\psi_{(g,h)}$ must vanish unless
\begin{gather*}
\varsigma_{(g,h)}\left(\sm{-1}00{-1}\right) = \xi_{g^{-1},h^{-1}}(k).
\end{gather*}
This is the analogue of the obstruction (2) in~\cite{Gaberdiel:2012gf}. We will employ this obstruction\footnote{In principle there is another independent possible obstruction coming from the projective class function property of $\psi_{(g,h)}$, namely $\psi_{(g,h)}$ must vanish unless $\xi_{(g,h)}(k)= 1$ for all $k$ commuting with both $g$ and $h$. However, in all cases apart from the $\ell=2$ case treated in \cite{Gaberdiel:2012gf}, this does not lead to new constraints.} in our computation and list the obstructed twisted-twined functions in Table~\ref{newfunctions}.

\looseness=-1 Note that in the case $\ell=2$, the first five conditions are equivalent to the ones discussed in~\cite{Gaberdiel:2012gf}.

\section{Groups, cohomologies, and representations}\label{Umbral Groups, Cohomologies, and Representations}

In this section we present the group-theoretic underpinnings of the theory. We introduce the umbral groups and discuss their cohomology, their rank-2 subgroups, and the projective representations of their centraliser subgroups.

\subsection{Group cohomology}\label{Group Cohomology}
The computation of the twisted-twined functions for generalised umbral moonshine depends crucially on the computation of the third group cohomology $H^3\big(\Gell,U(1)\big)$.

The cohomology groups are computed for all of the umbral groups using \textsc{GAP} and the module \textsc{HAP}. The results are listed in Table~\ref{groupstable}, where we include the root system $X$, the lambency $\ell$, the umbral group $\Gell$, of the corresponding cases of umbral moonshine. Throughout this paper we write $\ZZ/N\ZZ$ as~$N$. Also listed in the table is $n$, given by the group automorphism $\Gell \cong n.{\overline G}^{(\ell)}$, where ${\overline G}^{(\ell)}$ is the group of permutations of the irreducible components of the root systems of the corresponding Niemeier lattice induced by the lattice automorphism group. See~\cite{MUM} for more details. Note that $n=2$ or $n=1$ for all lambencies except for $\ell=6+3$ which has $n=3$.

The subgroups $n\subset \Gell$ play a special role in umbral moonshine: the corresponding twined functions are proportional to the untwined functions, $H_{(e,z),r} =c_r H_{(e,e),r} $ for a $c_r\in \CC$ for all $z\in n$ and $r\in \ZZ/2m\ZZ$. More generally, we have $H_{(e,zg),r} =c_r\,H_{(e,g),r} $ for all $g\in \Gell$. Together with the modularity condition~\eq{sl2z}, it suggests that the generalised umbral moonshine function~$H_{(g,h)}$ is really ``new" if $(g,h)$ cannot be obtained from the pair $(z,g')$ for some $z\in n$ via an~${\rm SL}_2(\ZZ)$ transformation. This inspires the following terminology: we say that a subgroup $\langle g, h \rangle$ is {\it old} if it is isomorphic to $\langle z, g' \rangle$ for some $g' \in G^{(\ell)}, z\in n$, and is {\it new} otherwise. We will also call the corresponding twisted-twined function an old resp. new function. In Table~\ref{groupstable} we also list $\lvert {\rm new} \rvert$, the number of conjugacy classes of rank-2 Abelian subgroups of $\Gell$ that are new.

The group theoretic computation for the case $\ell=2$ has been performed in \cite{Gaberdiel:2012gf}. For $\ell=3$, direct evaluation of cochains is computationally challenging. However, restriction of cochains to the $p$-Sylow subgroups $S_p(G)$ together with the inclusion map (cf.~\cite[Chapter~III-10]{MR1324339})
\begin{gather*}
H^3\big(\Gell,U(1)\big)\rightarrow \bigoplus_{p|\exp(G)} H^3\big(S_p\big(\Gell\big),U(1)\big),
\end{gather*}
allows us to simplify our calculations for $\ell=3$. In the above, $\exp(G)$ denotes the exponent of the group, and $\exp(G^{(3)})=1320$.

\begin{table}[t]\centering
\captionsetup{font=small}
\caption{Umbral groups and group cohomology.}\label{groupstable}
\newcolumntype{x}{>{\centering\arraybackslash\hspace{0pt}}p{0.8in}}
\resizebox{\columnwidth}{!}{
\begin{tabular}{c|xxxxxxxx}
$\rs$						&$A_1^{24}$		&$A_2^{12}$		&$A_3^8$			&$A_4^6$			&$A_5^4D_4$		&$A_6^4$			&$A_7^2D_5^2$	&$A_8^3$			\\
\vgap\doublerulesep\cline{1-9}\vgap\doublerulesep
$\ll$						&2				&3				&4				&5				&6				&7				&8				&9				\\
\vgap\doublerulesep\cline{1-9}\vgap\doublerulesep
$G^{\rs}$					&$M_{24}$		&$2.M_{12}$		&$2.\AGL_3(2)$	&$\GL_2(5)/2$		&$\GL_2(3)$		&$\SL_2(3)$		&$\Dih_4$			&$\Dih_6$			\\
$\bar{G}^{\rs}$				&$M_{24}$		&$M_{12}$		&$\AGL_3(2)$		&$\PGL_2(5)$		&$\PGL_2(3)$		&$\PSL_2(3)$		&$2^2$			&$\Sym_3$		\\
$H^3(\Gell,\mathbb{C}^\times)$&12				&$8\oplus 24$		&$4\oplus 24$		&$4\oplus 12$		&$2\oplus 24$		&24				&$2\oplus 2\oplus 4$	&$2\oplus 2\oplus 6$\\
$n$						&1				&2				&2				&2				&2				&2				&2				&2				\\
$\lvert$new$\rvert$ 			&34				&12				&8				&3				&0				&0				&0				&0				\\
\multicolumn{1}{r}{}\\
$\rs$						&$A_9^2D_6$		&$A_{11}D_7E_6$	&$A_{12}^2$		&$A_{15}D_9$		&$A_{17}E_7$		&$A_{24}$		&$D_4^{6}$		&$D_6^{4}$		\\
\vgap\doublerulesep\cline{1-9}\vgap\doublerulesep
$\ll$						&10				&12				&13				&16				&18				&25				&6+3			&10+5			\\
\vgap\doublerulesep\cline{1-9}\vgap\doublerulesep
$G^{\rs}$					&$4$				&$2$				&$4$				&$2$				&$2$				&$2$				&$3.\Sym_6$		&$\Sym_4$		\\
$\bar{G}^{\rs}$				&$2$				&$1$				&$2$				&$1$				&$1$				&$1$				&$\Sym_6$		&$\Sym_4$		\\
$H^3(\Gell,\mathbb{C}^\times)$&4				&1				&4				&1				&1				&1				&$2\oplus6\oplus 12$&$2\oplus12$		\\
$n$						&2				&2				&2				&2				&2				&2				&3				&2				\\
$\lvert$new$\rvert$ 			&0				&0				&0				&0				&0				&0				&8				&0				\\ 
\multicolumn{1}{r}{}\\
$\rs$						&$D_8^3$			&$D_{10}E_7^2$	&$D_{12}^2$		&$D_{16}E_8$		&$D_{24}$		&$E_6^4$			&$E_8^3$			\\
\vgap\doublerulesep\cline{1-8}\vgap\doublerulesep
$\ll$						&14+7			&18+9			&22+11			&30+15			&46+23			&12+4			&30+6,\,10,\,15 	\\
\vgap\doublerulesep\cline{1-8}\vgap\doublerulesep
$G^{\rs}$					&$\Sym_3$		&$2$				&$2$				&$1$				&$1$				&$\GL_2(3)$		&$\Sym_3$		\\
$\bar{G}^{\rs}$				&$\Sym_3$		&$2$				&$2$				&$1$				&$1$				&$\PGL_2(3)$		&$\Sym_3$		\\		
$H^3(\Gell,\mathbb{C}^\times)$&6				&1				&1				&1				&1				&$2\oplus24$		&6				\\
$n$						&1				&1				&1				&1				&1				&2				&1				\\
$\lvert$new$\rvert$ 			&0				&0				&0				&0				&0				&0				&0				\\ 
\end{tabular}}
\end{table}

\subsection{The 3-cocycles}

For most of the lambencies $\ell$ of umbral moonshine, the cohomology class of the \mbox{3-cocycle} $\o$ can be uniquely determined by the multiplier systems of the old functions $\psi^{}_{(e,g)}$ in the following way. Fix a lambency~$\ell$. For any $g\in G$, $\varsigma_{(e,g)}$ (cf.\ equation~\eqref{sl2z}) can be restricted to $\Gamma_0(\operatorname{ord}g)$ to obtain a multiplier system $ \varsigma_{(e,g)}\colon \Gamma_0(\operatorname{ord}g) \to \CC^\times$, which we assume to arise from a 3-cocycle $\o\colon G\times G\times G\to \CC^\times$ as discussed in Section~\ref{GUM}. At the same time, these multipliers are required to coincide with those of the known functions $\psi_{(e,g)}$~\cite{UM,MUM} (cf.\ Section~\ref{GUM}, condition~V). In this way, we obtain a consistency condition on~$[\o]$. This consistency condition uniquely fixes $[\o]\in H^3(G,U(1))$ for all umbral moonshine cases except for $\ell=3, 4, 7$. In these cases, there are exactly two cohomology classes that are consistent with the multiplier phases of the known functions $\psi_{(e,g)}^{}$ for all $g\in G$.

For the $\ell=4$ case, the two cohomology classes provide different multiplier systems on a~weight~1 Jacobi form associated to one of the new functions. One of the corresponding spaces of weight~1 weak Jacobi forms has dimension 0, but the vanishing of that specific new function would be incompatible with a decomposition of the characters into projective representations (cf.\ Section~\ref{GUM}, condition~IV). We therefore conclude that the new function must be non-zero and this allows us to uniquely select a 3-cohomology class.

For $\ell=3$ and $\ell=7$, both of the 3-cohomology classes that are consistent with the known functions $\psi_{(e,g)}^{}$ are also consistent with all the remaining twisted-twined functions we propose. As a result, in this sense generalised umbral moonshine for these two cases is compatible with both choices.

\subsection{The rank two subgroups}
Given a 3-cocycle $\o$, the projective class function property (cf.\ equation \eqref{proj_class_function}) guarantees that the generalised umbral moonshine function $\psi_{(g,h)}$ is determined by $\psi_{(g',h')}$ whenever the groups $\langle g,h\rangle$ and $\langle g',h'\rangle$ are conjugate to each other. It therefore suffices to specify the twisted-twined functions $\psi_{(g,h)}$ for a set of representatives in the coset
$\{g,h\in G\,\lvert\, gh=hg\}/{\sim}$, where the equivalence relation is generated by ${\rm SL}_2(\ZZ)$ transformations and by conjugation. In other words, we let $(g,h)\sim \gamma \big(g\lvert^k,h\lvert^k\big)$ for all $k\in G$ and $\g\in {\rm SL}_2(\ZZ)$.

 We tabulate these rank-2 subgroups in Table \ref{newfunctions} and indicate the new rank-2 subgroups whose twisted-twined functions are obstructed by the condition discussed at the end of Section~\ref{GUM}.
 In this table we also list the congruence subgroups $\Gamma_{(g,h)} \subseteq \text{SL}_2(\bz)$ stabilising the pair $(g,h)$ up to simultaneous conjugation (cf.~(\ref{def:gamma_gh})). In this list we employed the following notation for the congruence subgroups
\begin{gather*}
{\rm H}(p,q,r)=\left\{\gamma \in \text{SL}_2(\bz)\,\lvert\, \gamma =\sm1001\xmod\sm pqrp \right\},
\end{gather*}
which is a group if and only if $p|qr$. We also use the familiar convention for specific congruence subgroups:
\begin{gather*}
{\Gamma(N):= {\rm H}(N,N,N)},\qquad {\Gamma_1(N):= {\rm H}(N,1,N)} ,\qquad {\rm and} \qquad {\Gamma_0(N):= {\rm H}(1,1,N)}.
\end{gather*}
Finally, we have
\begin{gather*}
\Gamma_{(4A,4c)}=\left\{\sm abcd\in \Gamma(2)\,\lvert\, a+b+c\equiv 1\, \text{mod}\, 4 \right \}.
\end{gather*}
Here and everywhere else we denote the conjugacy class names of group elements $g\in G^{}$ with upper case letters, while the conjugacy class names of elements $h\in C_{G^{}}(g)$ in centralizer subgroups are denoted with lower case letters. See the file accompanying this article for an explicit construction of the umbral groups.

\subsection{Projective representations}

In this subsection we briefly explain how we construct the relevant projective representations. Our discussion follows closely Appendices~C and~D of~\cite{Gaberdiel:2012gf}.

Recall the definition of a projective representation in Section~\ref{Deformed Quantum Double and its Representations}. We say that two projective representations $\rho$ and $\rho'$ of a group $H$ are equivalent if they differ by a 1-cochain $\xi \in Z^1(\Gg,U(1))$, namely when $\rho(x) = \xi(x)\rho'(x)$ for all $x\in \Gg$. As a result, different classes in the second group cohomology $H^2(\Gg,U(1))$ lead to inequivalent projective representations. A convenient way to study the projective representations of a group $\Gg$ is to study the representations of a Schur cover $\Hh$ of $H$, defined as a central extension
\begin{gather*}
1\to M(\Gg)\to \Hh \xrightarrow{\phi} \Gg\to 1 ,
\end{gather*}
where the normal subgroup is the Schur multiplier $M(\Gg)\cong H_2(\Gg,\ZZ)$. As we will now explain, every representation of the Schur cover corresponds to a projective representation of $\Gg$ and vice versa.

Let $\phi^{-1}\colon \Gg\to \Hh$ be a lift of $H$. Let $\bar\rho\colon \Hh \to \operatorname{End}(V)$ be a representation of $\Hh$. Then $\rho:=\bar\rho\circ\phi^{-1}$ is a projective representation of $\Gg$ with the associated cocycle $c^\rho$ satisfying
\begin{gather*}
c^\rho(h_1,h_2) {{\bf 1}_V} = \overline{\rho} \left(\f^{-1}(h_1) \f^{-1}( h_2) \f^{-1}\big((h_1h_2)^{-1}\big)\right).
\end{gather*}
As a result, the 2-cocycle $c^\rho$ is fixed by $\bar\rho$ restricted to~${M(\Gg)}$. Every projective representation of~$\Gg$ is equivalent to a representation arising from its Schur cover in this way, and different (necessarily isoclinic) Schur covers of $\Gg$ lead to different but equivalent projective $\Gg$-representations.

Given a two-cohomology class $[c] \in H^2(\Gg,U(1))$, we would like to select the compatible representations from all the representations of $\Hh$. In other words, we would like to pick out the $\Hh$-representations that give rise to projective representations of $H$ whose corresponding 2-cocycle is consistent with $[c]$. Now we describe how this can be achieved. For a~representation~$\bar\rho$ of~$\Hh$ and $k\in M(\Gg)$, if $h_1,h_2\in \Gg$ satisfy $k=\phi^{-1}(h_1)\phi^{-1}(h_2)\phi^{-1}\big(h_1^{-1}\big)\phi^{-1}\big(h_2^{-1}\big)$, then $h_1$ and $h_2$ necessarily commute, and so $c^\rho(h_1,h_2) /c^\rho(h_2,h_1)$ depends only on the cohomology $[c^\rho]$. Moreover, the quotient can be evaluated as
\begin{gather}
 c^\rho(h_1,h_2) /c^\rho(h_2,h_1) {\bf 1}_V \nonumber\\
\qquad{} = \bar\rho\left(\phi^{-1}(h_1)\phi^{-1}(h_2)\phi^{-1}(h_1h_2)^{-1}\right) \bar\rho\left(\phi^{-1}(h_2)\phi^{-1}(h_1)\phi^{-1}(h_2h_1)^{-1}\right)^{-1}
 =\bar\rho(k),\label{filter_rep}
\end{gather}
which enables us to identify exactly the projective representations that are consistent with the given cohomology~$[c^\rho]$. Finally, irreducible representations of $\Hh$ give rise to irreducible projective representations of $\Gg$.

For our purpose, the finite groups whose projective representations we are interested in are the centralizer groups $C_G(g)$ where~$G$ is one of the umbral groups. The relevant 2-cocycle is determined by the 3-cocycle $\o \in Z^3(G,U(1))$ as in \eq{def:2cycle}. By considering the Schur cover of~$C_G(g)$ we obtain all irreducible projective representations~$C_G(g)$, and using \eq{filter_rep} these projective representations can then be filtered to those consistent with $[\o] \in H^3(G,U(1))$.

\section{The generalised umbral mock modular forms}\label{sec:mmf}

In this section, we employ the mock modularity property discussed in Section~\ref{GUM} to compute the mock modular forms of generalised umbral moonshine, which conjecturally arise from the generalised umbral moonshine module for the underlying deformed quantum double ${D^\o\big(\Gell\big)}$. In Section~\ref{Old Groups2} we compute the $g$-twisted $h$-twined functions when $\langle g,h\rangle$ is an old group. In Section~\ref{New Groups} we compute the remaining functions, the $g$-twisted $h$-twined functions when $\langle g,h\rangle$ is a new group.

\subsection{Old groups}\label{Old Groups2}

The untwisted twining functions of umbral moonshine were given in \cite{Cheng2010_1,UM,MUM,Eguchi2010a,Gaberdiel2010,Gaberdiel2010a}, and explicit expressions for them can be found in~\cite{proof}. For a given $\ell$ with index~$m$, these functions, which are denoted $H_{(e,g)}=(H_{(e,g),r})$, are $2m$-dimensional vector-valued mock modular forms of weight-1/2 which are related to the meromorphic Jacobi forms $\psi_{(e,g)}$ as discussed in Sections~\ref{UMreview} and~\ref{GUM}. In this subsection we will discuss how to obtain the meromorphic Jacobi form $\psi_{(g,h)}$ and the vector-valued mock modular form $H_{(g,h)}$ from the untwisted function $\psi_{(e,g')}$ for these cases.

 In the case there exists an ${\rm SL}_2(\ZZ)$ element $\g$ such that $\gamma(e,g') = (g,h)$, by the modularity condition \eq{sl2z}, we have
 \begin{gather}\label{sl2z_oldrankone}
\psi_{(g,h)}= \varsigma_{(e,g')}(\g) \psi_{(e,g')}|_{1,m}\gamma,
 \end{gather}
which can be evaluated explicitly from the relation between $\psi_{(e,g')}$ and $\psi_{(e,e)}$ and the explicit expression for $\psi_{(e,g')}$ given in~\cite{MUM,proof}.

When $\Gell \cong {\bar G}^{(\ell)}$ has $n=1$, all old groups are isomorphic to a cyclic group generated by an element $g'$, and the mock modular forms are found as described above. First consider umbral groups with $n=2$. Let $z\in n$ be the nontrivial element and note that $n$ is in the center of~$\Gell$. For the rank two old groups, namely those groups $\langle g,h\rangle$ that are isomorphic to $\langle z,g'\rangle$ for some $g'\in\Gell$, we provide a conjecture and a justification for the corresponding twisted-twined functions~$\psi_{(g,h)}$. From the modularity property \eq{sl2z}, we see that it suffices to calculate $\psi_{(z,g)}$. Fix the lambency $\ell$ with index $m$. We conjecture that
\begin{gather}\label{z_conj}\psi_{(z,g)} = \varsigma\psi_{(e,g)}|_{1,m}\big(\tfrac{1}{ 2}, 0\big)\end{gather}
for all $g\in G$ for some~$\varsigma\in\CC$.

To justify the conjecture, we start by verifying \eq{z_conj} whenever $\operatorname{ord} (g)$ is odd. Let $\gamma$ be an element of ${\rm SL}_2(\ZZ)$ such that $\gamma(e,zg)=(z,g)$. For instance, we can take $\g=\sm11{\ordg}{\ordg+1}$. Since $g$ is an element of the umbral group with odd order, the function $\psi_{(e,g)}$ transforms as a~weight one Jacobi form for $\Gamma_0(\ordg)$ with a multiplier system~\cite{MUM}. Hence we have
 \begin{gather*}
 \psi_{(z,g)} = \varsigma_{(e,zg)}(\g) \psi_{(e,zg)}\big\lvert_{1,m}\g = - \varsigma_{(e,zg)}(\g) \psi_{(e,g)}\big\lvert_m\big(0,\tfrac{1}{2}\big)\big\lvert_{1,m}\g
 =- {\varsigma_{(e,zg)}(\g)\over\varsigma_{(e,g)}(\g) }\psi_{(e,g)} \big\lvert_m\big(\tfrac{1}{2},0\big),
 \end{gather*}
 where we have used the operator relation \eq{modularelliptic} and the fact that
 \begin{gather}\label{shift_z2}
 \psi_{(e,zg)}= - \psi_{(e,g)}\big\lvert_m\big(0,\tfrac{1}{2}\big)
 \end{gather}
for all $g\in \Gell$.

Using $ \varsigma_{(z,g)}\left(\sm1101\right)\psi_{(z,g)}\big\lvert_{1,m}\sm1101 = \psi_{(z,zg)}$, it is easy to also verify the validity of \eq{z_conj} for elements $zg$ that are ``paired'' with the odd order elements $g$.
Finally, we conjecture that \eq{z_conj} is also true when $g$ and $zg$ both have even order.

This conjecture is strongly motivated by the finite group property discussed in Section~\ref{GUM}, which states that there exists an infinite dimensional module for $C_G(z) \cong G$ that underlies~$H_{(z,g)}$. Since replacing $\th_m$ by $\th_m\big\lvert_m\big(\tfrac{1}{2},0\big)$ simply reorders the components of the vector-valued mock modular forms $H_{(e,g)}=(H_{(e,g),r})$ by a shift $r\mapsto r+m$, the group theoretic consideration strongly suggests that this reordering has to happen for all $g\in G,$ given that it is true for all odd order elements $g$ as well as the elements $zg$ that they are paired with.

For $\ell=6+3$, the unique instance of umbral moonshine with $n=3$, the conjecture is slightly different. First, instead of \eq{shift_z2} we have now
\begin{gather*}
\psi_{(e,zg)} = \psi_{(e,g)}-\frac12\sum_{A\in \ZZ/3\ZZ} \psi_{(e,g)}|_{1,m}(0,A/3).
\end{gather*}
for $g\in\Gell$ with order $n_g$ not divisible by 3 that commutes with $z$. Choose any $\g = \sm abcd$ satisfying
\begin{gather*}
c= 0 \xmod{n_g} = 1 \xmod{3},\qquad d= 1\xmod{n_g} = 0 \xmod{3},
\end{gather*}
we have
\begin{gather*}
\psi_{(z,g)} = {\varsigma_{(e,zg)}(\g)\over\varsigma_{(e,g)}(\g) }\left( \psi_{(e,g)}-\frac12\sum_{A\in \ZZ/3\ZZ} \psi_{(e,g)}|_{1,m}(A/3,0)\right).
\end{gather*}

Repeating the same argument as before by evoking the projective representation property of moonshine, we conjecture that the above equality also holds for all $g$ commuting with~$z$, and hence for all rank two Abelian old groups $\langle z, g \rangle$.

\subsection{New groups}\label{New Groups}

Let $\langle g,h\rangle$ be a new group and consider the problem of evaluating $H_{(g,h)}=(H_{(g,h),r})$. First we discuss their analytic properties. We wish to show that for a new group $\langle g,h\rangle$, $H_{(g,h)}$ is bounded at all cusps.

Recall that the umbral moonshine conjecture \cite{UM,MUM} states that, for all $g\in G$, the untwisted twined functions $H_{(e,g)}$ are bounded at all cusps inequivalent to $\{{\rm i}\inf\}$ under the action of $\G_0(n_g)$. The number $n_g$ is given by the order of the group element $g$ as permutations of the irreducible components of the root system of the corresponding Niemeier lattice. See the appendices of~\cite{MUM} for the list of~$n_g$. In particular, this means that they have no poles at the cusp $\t\to 0$ for all $g\not\in n$ and the modular property of generalised umbral moonshine \eq{sl2z} implies that the corresponding twisted untwined functions $H_{(g,e)}$ are bounded at the cusp ${\rm i}\inf$. This fact together with the finite group property of generalised umbral moonshine~\eq{finite_grp_property} then implies that all twisted twined functions $H_{(g,h)}$ are bounded at the cups ${\rm i}\inf$ for all $g\in G$, $g\not\in n$ and for all $h\in C_G(g)$. Finally, evoking again the modular property \eq{sl2z} and considering $\g\in {\rm SL}_2(\ZZ)$ with $(g,h)= \gamma(g',h')$ to transform between cusps shows that the twisted twined functions $H_{(g,h)}$ are bounded at all cusps whenever $\langle g,h\rangle$ is a new group.

Next we discuss the modular properties of the new functions. From the above comment that~$H_{(g,h)}$ is bounded at all cusps and the absence of the mock Jacobi form in the corresponding space with non-trivial shadow (see the proof of Proposition~4.1.1 of~\cite{omjt}), we expect that the shadow of~$H_{(g,h)}$ vanishes. As a result, we conclude that $H_{(g,h)}$ is a weight $1/2$ vector-valued modular form for $\Gamma_{(g,h)}$ with the multiplier system determined by the 3-cocyle $\omega$ (cf.\ equation~\eqref{sl2z}). Note that any congruence subgroup, for instance $\G_{(g,h)}$, contains a subgroup of the form~$\Gamma(N)$. Hence each individual component $H_{(g,h),r}$ is a weight $1/2$ modular form (with trivial multiplier) for some $\Gamma(N)$. As a result, evoking Corollary~3 of~\cite{MR0472707}, we conclude that each component~$H_{(g,h),r}$ lies in the linear span of the theta functions
\begin{gather*}
\th^0_{m',s}(\t) :=\th_{m',s}(\t,0)=\sum_{k\in\ZZ} q^{\frac{(2m'k+s)^2}{4m'}}
\end{gather*}
for $m'\in \ZZ_+$, $s\in \ZZ/2m'\ZZ$.

\begin{table}[t!]\centering\small
\caption{New umbral subgroup representatives $\langle g,h\rangle\subseteq\Gell$.}\label{newfunctions}
 \newcommand{\gh}[2]{$\langle \text{#1},\text{#2}\rangle$}

\begin{tabular}{cllcl}
$\ell$ &$\langle g,h\rangle$& $\Gamma_{(g,h)}$ & Obs. 							& $\psi^{(\ell)}_{(g,h)}(\tau,\z)$   \\
\vgap\doublerulesep\cline{1-5}\vgap\doublerulesep
 \multirow{12}*{3}&$\langle \text{2B},\text{2a}\rangle$ & $\text{SL}_2(\bz)$ &  \checkmark	&0 \\
 &$\langle \text{2B},\text{2c}\rangle$ & $\Gamma_1(2)$ & \checkmark 		&0\\
 &$\langle \text{2B},\text{4a}\rangle$ & $\text{H}(2,2,4)$ & \checkmark 		&0\\
 &$\langle \text{2B},\text{4c}\rangle$ & $\text{H}(2,1,4)$ & \checkmark		&0\\
& $\langle \text{2C},\text{2b}\rangle$ & $\Gamma(2)$ & \checkmark			&0\\ 
& $\langle \text{2C},\text{2e}\rangle$ & $\text{SL}_2(\bz)$ & \checkmark		&0\\
& $\langle \text{3A},\text{3a}\rangle$ & $\text{SL}_2(\bz)$  &				&$-2  \theta^0_{4,2}\cdot\hat1+\big(\theta^0_{4,0}+\th^0_{4,4}\big)\cdot\hat{2}$\\
& $\langle \text{3A},\text{3b}\rangle$ & $\text{SL}_2(\bz)$ &				&$-2  \theta^0_{4,2}\cdot\hat1+\big(\theta^0_{4,0}+\th^0_{4,4}\big)\cdot\hat{2}$\\
& $\langle \text{3A},\text{3c}\rangle$ & $\Gamma_1(3)$ &					&$  \big(\theta^0_{9,3}-\theta^0_{9,9}\big)\cdot\hat1+\big(\theta^0_{9,0}-\theta^0_{9,6}\big)\cdot\hat2$\\
& $\langle \text{3A},\text{6a}\rangle$ & $\Gamma_1(6)$ &					&$   -\big(\theta^0_{9,3}-\theta^0_{9,9}\big)\cdot\hat1+\big(\theta^0_{9,0}-\theta^0_{9,6}\big)\cdot\hat2$\\
& $\langle \text{3A},\text{6b}\rangle$ & $\Gamma_1(2)$  &				& $2  \theta^0_{4,2}\cdot\hat1+\big(\theta^0_{4,0}+\th^0_{4,4}\big)\cdot\hat{2}$\\
& $\langle \text{3A},\text{6h}\rangle$& $\Gamma_1(2)$  &					& $2  \theta^0_{4,2}\cdot\hat1+\big(\theta^0_{4,0}+\th^0_{4,4}\big)\cdot\hat{2}$\\
\vgap\doublerulesep\cline{1-5}\vgap\doublerulesep
 \multirow{8}*{4}&  $\langle \text{2B},\text{2a}\rangle$  & $\Gamma_1(2)$ & \checkmark&0\\
&  $\langle \text{2B},\text{2b}\rangle$  & $\text{SL}_2(\bz)$  &\checkmark		&0\\ 
&  $\langle \text{2B},\text{4b}\rangle$  & $\Gamma_1(2)$  &				&$2\big(\theta^0_{4,0}-\theta^0_{4,4}\big)\cdot\hat2 $\\
&  $\langle \text{2C},\text{2a}\rangle$  & $\text{SL}_2(\bz)$  &\checkmark		&0\\
 & $\langle \text{2C},\text{2c}\rangle$  & $\text{SL}_2(\bz)$  &				&0\\
 & $\langle \text{2C},\text{2i}\rangle$   & $\Gamma_1(2)$ & \checkmark 		&0\\
 & $\langle \text{2C},\text{4c}\rangle$  & $\Gamma_1(2)$  &\checkmark		&0\\
 & $\langle \text{4A},\text{4c}\rangle$  & $\Gamma_{(4A,4c)}$ &			&$ (2-2{\rm i})  \theta^0_{2,1}\cdot\hat{2}$\\
\vgap\doublerulesep\cline{1-5}\vgap\doublerulesep
\multirow{3}*{5}&     $\langle \text{2B}, \text{2d}\rangle$ & $\text{SL}_2(\bz)$ & \checkmark&0  \\ 
&     $\langle \text{2B}, \text{2e}\rangle$ & $\Gamma_1(2)$  & \checkmark		&0\\
&    $\langle \text{2B}, \text{4a}\rangle$ & $\Gamma_1(4)$ &   				&$2 {\rm i} \big(\theta^0_{20,2}-\theta^0_{20,18}\big)\cdot\hat2+2{\rm i} \big(\theta^0_{20,6}-\theta^0_{20,14}\big)\cdot\hat4$\\
\vgap\doublerulesep\cline{1-5}\vgap\doublerulesep
\multirow{8}*{6+3}&    $\langle \text{2A}, \text{2a}\rangle$ & $\text{SL}_2(\bz)$ &\checkmark &0 \\
&    $\langle \text{2A}, \text{2b}\rangle$ & $\Gamma_1(2)$ &\checkmark &0 \\
&    $\langle \text{2A}, \text{2c}\rangle$ & $\Gamma_1(2)$ &\checkmark &0 \\
&    $\langle \text{2A}, \text{2b}\rangle$ & $\Gamma(2)$ &\checkmark &0 \\
&    $\langle \text{2A}, \text{6a}\rangle$ & $\Gamma_0(3)$ &\checkmark &0 \\
&    $\langle \text{2A}, \text{6a}\rangle$ & $\Gamma_0(3)$ &\checkmark &0 \\
&    $\langle \text{2A}, \text{2a}\rangle$ & $\text{SL}_2(\bz)$ &\checkmark &0 \\
&    $\langle \text{2B}, \text{4b}\rangle$ & $\text{H}(2,2,4)$ &\checkmark &0 \\
\vgap\doublerulesep\cline{1-5}
\end{tabular}

\end{table}

\looseness=1 Next we would like to further constrain the space of potential new functions of generalised umbral moonshine. Consider the congruence condition, which states that the set of $q$-powers of $H_{(g,h),r}$ with non-vanishing coefficients must be a subset of that of $H_{(g,e),r}$, a property that follows from the finite group property~\eq{finite_grp_property}. This defines a subspace $\Theta_{(g,h),r}$ generated by all theta functions $\th^0_{m',s}$ that $H_{(g,h),r}$ must lie in. Note that till now the space $\Theta_{(g,h),r}$ is independent of the twining element~$h$, as long as it generates a new group with the given twisting element~$g$. We can further restrict this space in two ways. First we require that also $H_{\gamma(g,h)}$ must have the correct $q$-powers for all $\g\in {\rm SL}_2(\ZZ)$. Second, elements of $\Gamma_{(g,h)}$ act as a linear transformation on this space, denoted $\sigma\colon \Gamma_{(g,h)}\to \operatorname{End}(\Theta_{(g,h),r})$, and $H_{(g,h),r}$ must lie in the simultaneous eigenspace of $\s(\g)$ for all $\g\in \Gamma_{(g,h)}$, with the eigenvalues given by the multiplier system $\xi_{(g,h)}$ (cf.\ equation~\eqref{proj_class_function}). In all of the unobstructed new groups, these simultaneous eigenspaces have dimension either~0 or~1. Together with the requirement of the existence of a~decomposition of the first homogeneous into projective representations, this fixes all the new functions.

In Table \ref{newfunctions} we also record all the new functions up to a phase. To avoid cluttering we have adopted the shorthand notation
\begin{gather*}
\hat i := \theta_{m,i}(\t,\z) - \theta_{m,-i}(\t,\z)
\end{gather*}
for $i\in \ZZ/2m\ZZ$.

\section{Summary and discussion}\label{Conclusions}

\looseness= 1 In this paper we introduced the notion of generalised umbral moonshine. Generalised umbral moonshine can be regarded as an extension of umbral moonshine and is analogous to the gene\-ra\-lised monstrous moonshine, albeit with significant conceptual and technical differences. When restricting to the~$A_1^{24}$ case it reduces to the generalised Mathieu moonshine analysed in~\cite{Gaberdiel:2012gf}.

We conjecture that the deformed Drinfel'd double of finite groups is the key structure underlying the (generalised) umbral moonshine phenomenon. More specifically, the generalised umbral moonshine conjecture states that for each of the 23 lambencies of umbral moonshine there is a way to assign an infinite-dimensional module for a deformed Drinfel'd double of the umbral group such that its graded character is given by the twisted-twined functions which we specified in Section~\ref{sec:mmf}. More explicitly, we proposed in Section~\ref{GUM} six conditions of generalised umbral moonshine conjecture and provided evidence in the form of proposing all the twisted-twined functions for all cases of generalised umbral moonshine. A~proof of this conjecture should be attaianable following the methods of~\cite{proof,MR3582425}.

\looseness = 1 It is interesting that one does not obtain new mock modular forms (with non-vanishing shadow) in this generalisation of umbral moonshine. As we reviewed in Section~\ref{UMreview}, the origin of the mockness of the umbral moonshine mock modular forms $H_{g}^{(\ell)}$ is the poles of the meromorphic Jacobi forms $\psi_{g}^{(\ell)}$ \eq{polarpart}, and the residues of the poles of these functions are given by characters of the umbral group. From \eq{sl2z} and the conjecture in Section~\ref{Old Groups2}, we see that the shadow of the mock modular form $H_{(g,h)}^{(\ell)}$ when $\langle g,h \rangle$ is an old group can be obtained from that of certain umbral moonshine function $H_{g'}^{(\ell)}$. In the case when $\langle g,h \rangle$ is a new group, the shadow of $H_{(g,h)}^{(\ell)}$ vanishes as a result of the reasoning presented in Section~\ref{New Groups}. In fact, all the new functions, which do not play a role in umbral moonshine itself, are given by theta functions (see Table~\ref{newfunctions}).

\looseness = 1 The existence of generalised umbral moonshine is particularly interesting in view of the absence (so far) of a clear relation between umbral moonshine and conventional CFTs\footnote{From this point of view, the fact that there is a generalised moonshine for each of the~23 cases of umbral moonshine is perhaps more surprising than the existence of generalised Mathieu moonshine alone. This is because, unlike the other 22 cases of umbral moonshine, Mathieu moonshine has an obvious relation to CFTs. Indeed, in~\cite{Gaberdiel:2012gf} the existence of generalised Mathieu moonshine is argued to suggest that a holomorphic VOA may be underlying Mathieu moonshine.}. As mentioned in Section~\ref{GMM}, the generalised monstrous moonshine is believed to be a manifestation of the fact that one can define the twisted modules in view of the existence of the moonshine VOA. More generally, the physical context of the deformed Drinfel'd double of a finite group is usually assumed to be that of holomorphic orbifold CFTs or closely related structures such as 3d TQFTs. The important lesson generalised umbral moonshine has taught us is the following: whatever algebraic structure underlies umbral moonshine, it has to retain the intimate relation to the deformed Drinfel'd double as chiral conformal field theories do.

\appendix
\section{Examples}

\looseness = 1 In this appendix we focus on three examples of $g$-twisted generalised umbral moonshine functions for lambency $\ell=3,4,6+3$ of index $m=3,4,6$ respectively. In particular, we list the independent components of the vector-valued mock modular forms $H_{(g,h)}=(H_{(g,h),r})$, $r\in \ZZ/2m\ZZ$, for all $h\in C_{G^{\ell}}(g)$. Our notation for the conjugacy classes can be found in the attached file. In these tables, the rational numbers in the first column denotes the $q$-power of the corresponding terms in the Fourier expansion of $H_{(g,h),r}$. We also list the projective character table for the projective representations compatible with the 3-cocycle, as well as the decompositions of the first few homogeneous components $K^g_{r,\a}$ of the conjectural generalised umbral moonshine module for these cases (cf.\ equation~\eqref{finite_grp_property}). To increase readability we use~``.'' to denote zero.

\newpage

\subsection[$\ell=3$, 4B twist]{$\boldsymbol{\ell=3}$, 4B twist}
\begin{table}[h!]\centering
\caption{Projective character table for $C_{G^{(3)}}(g)$.}
\resizebox{\columnwidth}{!}{
\begin{tabular}{r|rrrrrrrrrrrrrrrrrrrr}
$$&$1a$&$2a$&$2b$&$2c$&$2d$&$2e$&$4a$&$4b$&$4c$&$4d$&$4e$&$4f$&$8a$&$8b$&$8c$&$8d$&$8e$&$8f$&$8g$&$8h$\\
\hline$\chi_{1}$&$1$&$1$&$1$&$1$&$1$&$1$&$1$&$1$&$1$&$1$&$1$&$1$&$1$&$1$&$1$&$1$&$1$&$1$&$1$&$1$\\
$\chi_{2}$&$1$&$-1$&$1$&$-1$&$1$&$1$&$1$&$-1$&$1$&$-1$&$1$&$1$&$-1$&$1$&$-1$&$1$&$-1$&$1$&$-1$&$1$\\
$\chi_{3}$&$1$&$1$&$1$&$1$&$1$&$1$&$1$&$1$&$1$&$1$&$1$&$1$&$-1$&$-1$&$-1$&$-1$&$-1$&$-1$&$-1$&$-1$\\
$\chi_{4}$&$1$&$-1$&$1$&$-1$&$1$&$1$&$1$&$-1$&$1$&$-1$&$1$&$1$&$1$&$-1$&$1$&$-1$&$1$&$-1$&$1$&$-1$\\
$\chi_{5}$&$1$&$-1$&$1$&$1$&$-1$&$-1$&$i$&$-i$&$-i$&$i$&$-i$&$i$&$\zeta_{8}^{7}$&$\zeta_{8}^{3}$&$\zeta_{8}^{5}$&$\zeta_{8}$&$\zeta_{8}^{3}$&$\zeta_{8}^{7}$&$\zeta_{8}$&$\zeta_{8}^{5}$\\
$\chi_{6}$&$1$&$-1$&$1$&$1$&$-1$&$-1$&$i$&$-i$&$-i$&$i$&$-i$&$i$&$\zeta_{8}^{3}$&$\zeta_{8}^{7}$&$\zeta_{8}$&$\zeta_{8}^{5}$&$\zeta_{8}^{7}$&$\zeta_{8}^{3}$&$\zeta_{8}^{5}$&$\zeta_{8}$\\
$\chi_{7}$&$1$&$-1$&$1$&$1$&$-1$&$-1$&$-i$&$i$&$i$&$-i$&$i$&$-i$&$\zeta_{8}^{5}$&$\zeta_{8}$&$\zeta_{8}^{7}$&$\zeta_{8}^{3}$&$\zeta_{8}$&$\zeta_{8}^{5}$&$\zeta_{8}^{3}$&$\zeta_{8}^{7}$\\
$\chi_{8}$&$1$&$-1$&$1$&$1$&$-1$&$-1$&$-i$&$i$&$i$&$-i$&$i$&$-i$&$\zeta_{8}$&$\zeta_{8}^{5}$&$\zeta_{8}^{3}$&$\zeta_{8}^{7}$&$\zeta_{8}^{5}$&$\zeta_{8}$&$\zeta_{8}^{7}$&$\zeta_{8}^{3}$\\
$\chi_{9}$&$1$&$1$&$1$&$-1$&$-1$&$-1$&$i$&$i$&$-i$&$-i$&$-i$&$i$&$\zeta_{8}^{7}$&$\zeta_{8}^{7}$&$\zeta_{8}^{5}$&$\zeta_{8}^{5}$&$\zeta_{8}^{3}$&$\zeta_{8}^{3}$&$\zeta_{8}$&$\zeta_{8}$\\
$\chi_{10}$&$1$&$1$&$1$&$-1$&$-1$&$-1$&$i$&$i$&$-i$&$-i$&$-i$&$i$&$\zeta_{8}^{3}$&$\zeta_{8}^{3}$&$\zeta_{8}$&$\zeta_{8}$&$\zeta_{8}^{7}$&$\zeta_{8}^{7}$&$\zeta_{8}^{5}$&$\zeta_{8}^{5}$\\
$\chi_{11}$&$1$&$1$&$1$&$-1$&$-1$&$-1$&$-i$&$-i$&$i$&$i$&$i$&$-i$&$\zeta_{8}^{5}$&$\zeta_{8}^{5}$&$\zeta_{8}^{7}$&$\zeta_{8}^{7}$&$\zeta_{8}$&$\zeta_{8}$&$\zeta_{8}^{3}$&$\zeta_{8}^{3}$\\
$\chi_{12}$&$1$&$1$&$1$&$-1$&$-1$&$-1$&$-i$&$-i$&$i$&$i$&$i$&$-i$&$\zeta_{8}$&$\zeta_{8}$&$\zeta_{8}^{3}$&$\zeta_{8}^{3}$&$\zeta_{8}^{5}$&$\zeta_{8}^{5}$&$\zeta_{8}^{7}$&$\zeta_{8}^{7}$\\
$\chi_{13}$&$1$&$-1$&$1$&$-1$&$1$&$1$&$-1$&$1$&$-1$&$1$&$-1$&$-1$&$-i$&$i$&$i$&$-i$&$-i$&$i$&$i$&$-i$\\
$\chi_{14}$&$1$&$-1$&$1$&$-1$&$1$&$1$&$-1$&$1$&$-1$&$1$&$-1$&$-1$&$i$&$-i$&$-i$&$i$&$i$&$-i$&$-i$&$i$\\
$\chi_{15}$&$1$&$1$&$1$&$1$&$1$&$1$&$-1$&$-1$&$-1$&$-1$&$-1$&$-1$&$-i$&$-i$&$i$&$i$&$-i$&$-i$&$i$&$i$\\
$\chi_{16}$&$1$&$1$&$1$&$1$&$1$&$1$&$-1$&$-1$&$-1$&$-1$&$-1$&$-1$&$i$&$i$&$-i$&$-i$&$i$&$i$&$-i$&$-i$\\
$\chi_{17}$&$2$&$.$&$-2$&$.$&$-2$&$2$&$2$&$.$&$2$&$.$&$-2$&$-2$&$.$&$.$&$.$&$.$&$.$&$.$&$.$&$.$\\
$\chi_{18}$&$2$&$.$&$-2$&$.$&$-2$&$2$&$-2$&$.$&$-2$&$.$&$2$&$2$&$.$&$.$&$.$&$.$&$.$&$.$&$.$&$.$\\
$\chi_{19}$&$2$&$.$&$-2$&$.$&$2$&$-2$&$2 i$&$.$&$-2 i$&$.$&$2 i$&$-2 i$&$.$&$.$&$.$&$.$&$.$&$.$&$.$&$.$\\
$\chi_{20}$&$2$&$.$&$-2$&$.$&$2$&$-2$&$-2 i$&$.$&$2 i$&$.$&$-2 i$&$2 i$&$.$&$.$&$.$&$.$&$.$&$.$&$.$&$.$\\
\hline
\end{tabular}}
\end{table}

\begin{table}[h!]\centering
\caption{Coefficients of $H_{(g,h),1}$.}
\resizebox{\columnwidth}{!}{
\begin{tabular}{r|rrrrrrrrrrrrrrrrrrrr}
$$&$1a$&$2a$&$2b$&$2c$&$2d$&$2e$&$4a$&$4b$&$4c$&$4d$&$4e$&$4f$&$8a$&$8b$&$8c$&$8d$&$8e$&$8f$&$8g$&$8h$\\
\hline$2/12$&$2$&$.$&$2$&$.$&$2$&$2$&$2$&$.$&$2$&$.$&$2$&$2$&$.$&$2$&$.$&$2$&$.$&$2$&$.$&$2$\\
$5/12$&$4$&$.$&$-4$&$.$&$4$&$-4$&$4 i$&$.$&$-4 i$&$.$&$4 i$&$-4 i$&$.$&$.$&$.$&$.$&$.$&$.$&$.$&$.$\\
$8/12$&$6$&$.$&$6$&$.$&$6$&$6$&$-6$&$.$&$-6$&$.$&$-6$&$-6$&$.$&$2 i$&$.$&$-2 i$&$.$&$2 i$&$.$&$-2 i$\\
$11/12$&$8$&$.$&$-8$&$.$&$8$&$-8$&$-8 i$&$.$&$8 i$&$.$&$-8 i$&$8 i$&$.$&$.$&$.$&$.$&$.$&$.$&$.$&$.$\\
$14/12$&$12$&$.$&$12$&$.$&$12$&$12$&$12$&$.$&$12$&$.$&$12$&$12$&$.$&$-4$&$.$&$-4$&$.$&$-4$&$.$&$-4$\\
$17/12$&$16$&$.$&$-16$&$.$&$16$&$-16$&$16 i$&$.$&$-16 i$&$.$&$16 i$&$-16 i$&$.$&$.$&$.$&$.$&$.$&$.$&$.$&$.$\\
$20/12$&$20$&$.$&$20$&$.$&$20$&$20$&$-20$&$.$&$-20$&$.$&$-20$&$-20$&$.$&$-4 i$&$.$&$4 i$&$.$&$-4 i$&$.$&$4 i$\\
$23/12$&$28$&$.$&$-28$&$.$&$28$&$-28$&$-28 i$&$.$&$28 i$&$.$&$-28 i$&$28 i$&$.$&$.$&$.$&$.$&$.$&$.$&$.$&$.$\\
$26/12$&$36$&$.$&$36$&$.$&$36$&$36$&$36$&$.$&$36$&$.$&$36$&$36$&$.$&$4$&$.$&$4$&$.$&$4$&$.$&$4$\\
$29/12$&$44$&$.$&$-44$&$.$&$44$&$-44$&$44 i$&$.$&$-44 i$&$.$&$44 i$&$-44 i$&$.$&$.$&$.$&$.$&$.$&$.$&$.$&$.$\\
\hline
\end{tabular}}
\end{table}

\begin{table}[h!]\centering
\caption{Coefficients of $H_{(g,h),2}$.}
\resizebox{\columnwidth}{!}{
\begin{tabular}{r|rrrrrrrrrrrrrrrrrrrr}
$$&$1a$&$2a$&$2b$&$2c$&$2d$&$2e$&$4a$&$4b$&$4c$&$4d$&$4e$&$4f$&$8a$&$8b$&$8c$&$8d$&$8e$&$8f$&$8g$&$8h$\\
\hline$2/12$&$2$&$.$&$2$&$.$&$-2$&$-2$&$2 i$&$.$&$-2 i$&$.$&$-2 i$&$2 i$&$2\zeta_{8}^{7}$&$.$&$2\zeta_{8}^{5}$&$.$&$2\zeta_{8}^{3}$&$.$&$2\zeta_{8}$&$.$\\
$5/12$&$4$&$.$&$-4$&$.$&$-4$&$4$&$-4$&$.$&$-4$&$.$&$4$&$4$&$.$&$.$&$.$&$.$&$.$&$.$&$.$&$.$\\
$8/12$&$6$&$.$&$6$&$.$&$-6$&$-6$&$-6 i$&$.$&$6 i$&$.$&$6 i$&$-6 i$&$2\zeta_{8}^{5}$&$.$&$2\zeta_{8}^{7}$&$.$&$2\zeta_{8}$&$.$&$2\zeta_{8}^{3}$&$.$\\
$11/12$&$8$&$.$&$-8$&$.$&$-8$&$8$&$8$&$.$&$8$&$.$&$-8$&$-8$&$.$&$.$&$.$&$.$&$.$&$.$&$.$&$.$\\
$14/12$&$12$&$.$&$12$&$.$&$-12$&$-12$&$12 i$&$.$&$-12 i$&$.$&$-12 i$&$12 i$&$4\zeta_{8}^{3}$&$.$&$4\zeta_{8}$&$.$&$4\zeta_{8}^{7}$&$.$&$4\zeta_{8}^{5}$&$.$\\
$17/12$&$16$&$.$&$-16$&$.$&$-16$&$16$&$-16$&$.$&$-16$&$.$&$16$&$16$&$.$&$.$&$.$&$.$&$.$&$.$&$.$&$.$\\
$20/12$&$20$&$.$&$20$&$.$&$-20$&$-20$&$-20 i$&$.$&$20 i$&$.$&$20 i$&$-20 i$&$4\zeta_{8}$&$.$&$4\zeta_{8}^{3}$&$.$&$4\zeta_{8}^{5}$&$.$&$4\zeta_{8}^{7}$&$.$\\
$23/12$&$28$&$.$&$-28$&$.$&$-28$&$28$&$28$&$.$&$28$&$.$&$-28$&$-28$&$.$&$.$&$.$&$.$&$.$&$.$&$.$&$.$\\
$26/12$&$36$&$.$&$36$&$.$&$-36$&$-36$&$36 i$&$.$&$-36 i$&$.$&$-36 i$&$36 i$&$4\zeta_{8}^{7}$&$.$&$4\zeta_{8}^{5}$&$.$&$4\zeta_{8}^{3}$&$.$&$4\zeta_{8}$&$.$\\
$29/12$&$44$&$.$&$-44$&$.$&$-44$&$44$&$-44$&$.$&$-44$&$.$&$44$&$44$&$.$&$.$&$.$&$.$&$.$&$.$&$.$&$.$\\
\hline
\end{tabular}}
\end{table}

\begin{table}[h!]\centering
\caption{Decompositions for the first few $K^g_{1,\a}.$}
\resizebox{\columnwidth}{!}{
\begin{tabular}{r|rrrrrrrrrrrrrrrrrrrr}
$$&$\chi_{1}$&$\chi_{2}$&$\chi_{3}$&$\chi_{4}$&$\chi_{5}$&$\chi_{6}$&$\chi_{7}$&$\chi_{8}$&$\chi_{9}$&$\chi_{10}$&$\chi_{11}$&$\chi_{12}$&$\chi_{13}$&$\chi_{14}$&$\chi_{15}$&$\chi_{16}$&$\chi_{17}$&$\chi_{18}$&$\chi_{19}$&$\chi_{20}$\\
\hline$2/12$&$1$&$1$&$.$&$.$&$.$&$.$&$.$&$.$&$.$&$.$&$.$&$.$&$.$&$.$&$.$&$.$&$.$&$.$&$.$&$.$\\
$5/12$&$.$&$.$&$.$&$.$&$.$&$.$&$.$&$.$&$.$&$.$&$.$&$.$&$.$&$.$&$.$&$.$&$.$&$.$&$2$&$.$\\
$8/12$&$.$&$.$&$.$&$.$&$.$&$.$&$.$&$.$&$.$&$.$&$.$&$.$&$2$&$1$&$1$&$2$&$.$&$.$&$.$&$.$\\
$11/12$&$.$&$.$&$.$&$.$&$.$&$.$&$.$&$.$&$.$&$.$&$.$&$.$&$.$&$.$&$.$&$.$&$.$&$.$&$.$&$4$\\
$14/12$&$2$&$2$&$4$&$4$&$.$&$.$&$.$&$.$&$.$&$.$&$.$&$.$&$.$&$.$&$.$&$.$&$.$&$.$&$.$&$.$\\
$17/12$&$.$&$.$&$.$&$.$&$.$&$.$&$.$&$.$&$.$&$.$&$.$&$.$&$.$&$.$&$.$&$.$&$.$&$.$&$8$&$.$\\
$20/12$&$.$&$.$&$.$&$.$&$.$&$.$&$.$&$.$&$.$&$.$&$.$&$.$&$4$&$6$&$6$&$4$&$.$&$.$&$.$&$.$\\
$23/12$&$.$&$.$&$.$&$.$&$.$&$.$&$.$&$.$&$.$&$.$&$.$&$.$&$.$&$.$&$.$&$.$&$.$&$.$&$.$&$14$\\
$26/12$&$10$&$10$&$8$&$8$&$.$&$.$&$.$&$.$&$.$&$.$&$.$&$.$&$.$&$.$&$.$&$.$&$.$&$.$&$.$&$.$\\
$29/12$&$.$&$.$&$.$&$.$&$.$&$.$&$.$&$.$&$.$&$.$&$.$&$.$&$.$&$.$&$.$&$.$&$.$&$.$&$22$&$.$\\
\hline
\end{tabular}}
\end{table}

\begin{table}[h!]\centering
\caption{Decompositions for the first few $K^g_{2,\a}$.}
\resizebox{\columnwidth}{!}{
\begin{tabular}{r|rrrrrrrrrrrrrrrrrrrr}
$$&$\chi_{1}$&$\chi_{2}$&$\chi_{3}$&$\chi_{4}$&$\chi_{5}$&$\chi_{6}$&$\chi_{7}$&$\chi_{8}$&$\chi_{9}$&$\chi_{10}$&$\chi_{11}$&$\chi_{12}$&$\chi_{13}$&$\chi_{14}$&$\chi_{15}$&$\chi_{16}$&$\chi_{17}$&$\chi_{18}$&$\chi_{19}$&$\chi_{20}$\\
\hline$2/12$&$.$&$.$&$.$&$.$&$1$&$.$&$.$&$.$&$1$&$.$&$.$&$.$&$.$&$.$&$.$&$.$&$.$&$.$&$.$&$.$\\
$5/12$&$.$&$.$&$.$&$.$&$.$&$.$&$.$&$.$&$.$&$.$&$.$&$.$&$.$&$.$&$.$&$.$&$.$&$2$&$.$&$.$\\
$8/12$&$.$&$.$&$.$&$.$&$.$&$.$&$2$&$1$&$.$&$.$&$2$&$1$&$.$&$.$&$.$&$.$&$.$&$.$&$.$&$.$\\
$11/12$&$.$&$.$&$.$&$.$&$.$&$.$&$.$&$.$&$.$&$.$&$.$&$.$&$.$&$.$&$.$&$.$&$4$&$.$&$.$&$.$\\
$14/12$&$.$&$.$&$.$&$.$&$2$&$4$&$.$&$.$&$2$&$4$&$.$&$.$&$.$&$.$&$.$&$.$&$.$&$.$&$.$&$.$\\
$17/12$&$.$&$.$&$.$&$.$&$.$&$.$&$.$&$.$&$.$&$.$&$.$&$.$&$.$&$.$&$.$&$.$&$.$&$8$&$.$&$.$\\
$20/12$&$.$&$.$&$.$&$.$&$.$&$.$&$4$&$6$&$.$&$.$&$4$&$6$&$.$&$.$&$.$&$.$&$.$&$.$&$.$&$.$\\
$23/12$&$.$&$.$&$.$&$.$&$.$&$.$&$.$&$.$&$.$&$.$&$.$&$.$&$.$&$.$&$.$&$.$&$14$&$.$&$.$&$.$\\
$26/12$&$.$&$.$&$.$&$.$&$10$&$8$&$.$&$.$&$10$&$8$&$.$&$.$&$.$&$.$&$.$&$.$&$.$&$.$&$.$&$.$\\
$29/12$&$.$&$.$&$.$&$.$&$.$&$.$&$.$&$.$&$.$&$.$&$.$&$.$&$.$&$.$&$.$&$.$&$.$&$22$&$.$&$.$\\
\hline
\end{tabular}}
\end{table}

\newpage

\subsection[$\ell=4$, 2B twist]{$\boldsymbol{\ell=4}$, 2B twist}
\begin{table}[h!]
\centering
\caption{Projective character table for $C_{G^{(4)}}(g)$.}
\resizebox{\columnwidth}{!}{
\begin{tabular}{r|rrrrrrrrrrrrrrrrr}
$$&$1a$&$2a$&$2b$&$2c$&$2d$&$2e$&$3a$&$3b$&$4a$&$4b$&$4c$&$6a$&$6b$&$6c$&$6d$&$6e$&$6f$\\
\hline$\chi_{1}$&$2$&$.$&$.$&$2$&$2$&$2$&$1$&$1$&$2$&$.$&$.$&$1$&$1$&$1$&$1$&$1$&$1$\\
$\chi_{2}$&$2$&$.$&$.$&$-2$&$2$&$-2$&$1$&$1$&$.$&$.$&$2$&$1$&$-1$&$-1$&$-1$&$-1$&$1$\\
$\chi_{3}$&$2$&$.$&$.$&$-2$&$-2$&$2$&$1$&$1$&$.$&$2$&$.$&$-1$&$-1$&$-1$&$1$&$1$&$-1$\\
$\chi_{4}$&$2$&$.$&$.$&$2$&$2$&$2$&$\zeta_{3}^{2}$&$\zeta_{3}$&$2$&$.$&$.$&$\zeta_{3}$&$\zeta_{3}^{2}$&$\zeta_{3}$&$\zeta_{3}^{2}$&$\zeta_{3}$&$\zeta_{3}^{2}$\\
$\chi_{5}$&$2$&$.$&$.$&$2$&$2$&$2$&$\zeta_{3}$&$\zeta_{3}^{2}$&$2$&$.$&$.$&$\zeta_{3}^{2}$&$\zeta_{3}$&$\zeta_{3}^{2}$&$\zeta_{3}$&$\zeta_{3}^{2}$&$\zeta_{3}$\\
$\chi_{6}$&$2$&$.$&$.$&$-2$&$2$&$-2$&$\zeta_{3}^{2}$&$\zeta_{3}$&$.$&$.$&$2$&$\zeta_{3}$&$\zeta_{6}$&$\zeta_{6}^{5}$&$\zeta_{6}$&$\zeta_{6}^{5}$&$\zeta_{3}^{2}$\\
$\chi_{7}$&$2$&$.$&$.$&$-2$&$2$&$-2$&$\zeta_{3}$&$\zeta_{3}^{2}$&$.$&$.$&$2$&$\zeta_{3}^{2}$&$\zeta_{6}^{5}$&$\zeta_{6}$&$\zeta_{6}^{5}$&$\zeta_{6}$&$\zeta_{3}$\\
$\chi_{8}$&$2$&$.$&$.$&$-2$&$-2$&$2$&$\zeta_{3}^{2}$&$\zeta_{3}$&$.$&$2$&$.$&$\zeta_{6}^{5}$&$\zeta_{6}$&$\zeta_{6}^{5}$&$\zeta_{3}^{2}$&$\zeta_{3}$&$\zeta_{6}$\\
$\chi_{9}$&$2$&$.$&$.$&$-2$&$-2$&$2$&$\zeta_{3}$&$\zeta_{3}^{2}$&$.$&$2$&$.$&$\zeta_{6}$&$\zeta_{6}^{5}$&$\zeta_{6}$&$\zeta_{3}$&$\zeta_{3}^{2}$&$\zeta_{6}^{5}$\\
$\chi_{10}$&$4$&$.$&$.$&$4$&$-4$&$-4$&$-1$&$-1$&$.$&$.$&$.$&$1$&$-1$&$-1$&$1$&$1$&$1$\\
$\chi_{11}$&$4$&$.$&$.$&$4$&$-4$&$-4$&$\zeta_{6}$&$\zeta_{6}^{5}$&$.$&$.$&$.$&$\zeta_{3}$&$\zeta_{6}$&$\zeta_{6}^{5}$&$\zeta_{3}^{2}$&$\zeta_{3}$&$\zeta_{3}^{2}$\\
$\chi_{12}$&$4$&$.$&$.$&$4$&$-4$&$-4$&$\zeta_{6}^{5}$&$\zeta_{6}$&$.$&$.$&$.$&$\zeta_{3}^{2}$&$\zeta_{6}^{5}$&$\zeta_{6}$&$\zeta_{3}$&$\zeta_{3}^{2}$&$\zeta_{3}$\\
$\chi_{13}$&$6$&$.$&$.$&$6$&$6$&$6$&$.$&$.$&$-2$&$.$&$.$&$.$&$.$&$.$&$.$&$.$&$.$\\
$\chi_{14}$&$6$&$.$&$.$&$-6$&$6$&$-6$&$.$&$.$&$.$&$.$&$-2$&$.$&$.$&$.$&$.$&$.$&$.$\\
$\chi_{15}$&$6$&$.$&$.$&$-6$&$-6$&$6$&$.$&$.$&$.$&$-2$&$.$&$.$&$.$&$.$&$.$&$.$&$.$\\
\hline
\end{tabular}}
\end{table}

\newpage

\begin{table}[h!]\centering
\caption{Coefficients of $H_{(g,h),1}$.}
\resizebox{\columnwidth}{!}{
\begin{tabular}{r|rrrrrrrrrrrrrrrrr}
$$&$1a$&$2a$&$2b$&$2c$&$2d$&$2e$&$3a$&$3b$&$4a$&$4b$&$4c$&$6a$&$6b$&$6c$&$6d$&$6e$&$6f$\\
\hline$3/16$&$4$&$.$&$.$&$4$&$4$&$4$&$-1$&$-1$&$4$&$.$&$.$&$-1$&$-1$&$-1$&$-1$&$-1$&$-1$\\
$11/16$&$4$&$.$&$.$&$-4$&$4$&$-4$&$2$&$2$&$.$&$.$&$4$&$2$&$-2$&$-2$&$-2$&$-2$&$2$\\
$19/16$&$12$&$.$&$.$&$12$&$12$&$12$&$.$&$.$&$-4$&$.$&$.$&$.$&$.$&$.$&$.$&$.$&$.$\\
$27/16$&$16$&$.$&$.$&$-16$&$16$&$-16$&$-1$&$-1$&$.$&$.$&$.$&$-1$&$1$&$1$&$1$&$1$&$-1$\\
$35/16$&$24$&$.$&$.$&$24$&$24$&$24$&$.$&$.$&$8$&$.$&$.$&$.$&$.$&$.$&$.$&$.$&$.$\\
$43/16$&$36$&$.$&$.$&$-36$&$36$&$-36$&$.$&$.$&$.$&$.$&$4$&$.$&$.$&$.$&$.$&$.$&$.$\\
$51/16$&$56$&$.$&$.$&$56$&$56$&$56$&$-2$&$-2$&$-8$&$.$&$.$&$-2$&$-2$&$-2$&$-2$&$-2$&$-2$\\
$59/16$&$76$&$.$&$.$&$-76$&$76$&$-76$&$2$&$2$&$.$&$.$&$-4$&$2$&$-2$&$-2$&$-2$&$-2$&$2$\\
$67/16$&$108$&$.$&$.$&$108$&$108$&$108$&$.$&$.$&$12$&$.$&$.$&$.$&$.$&$.$&$.$&$.$&$.$\\
$75/16$&$148$&$.$&$.$&$-148$&$148$&$-148$&$-1$&$-1$&$.$&$.$&$4$&$-1$&$1$&$1$&$1$&$1$&$-1$\\
\hline
\end{tabular}}
\end{table}

\begin{table}[h]
\centering
\caption{Coefficients of $H_{(g,h),2}$.}
\resizebox{\columnwidth}{!}{
\begin{tabular}{r|rrrrrrrrrrrrrrrrr}
$$&$1a$&$2a$&$2b$&$2c$&$2d$&$2e$&$3a$&$3b$&$4a$&$4b$&$4c$&$6a$&$6b$&$6c$&$6d$&$6e$&$6f$\\
\hline$0/2$&$2$&$.$&$.$&$-2$&$-2$&$2$&$1$&$1$&$.$&$2$&$.$&$-1$&$-1$&$-1$&$1$&$1$&$-1$\\
$1/2$&$8$&$.$&$.$&$8$&$-8$&$-8$&$-2$&$-2$&$.$&$.$&$.$&$2$&$-2$&$-2$&$2$&$2$&$2$\\
$2/2$&$12$&$.$&$.$&$-12$&$-12$&$12$&$.$&$.$&$.$&$-4$&$.$&$.$&$.$&$.$&$.$&$.$&$.$\\
$3/2$&$16$&$.$&$.$&$16$&$-16$&$-16$&$2$&$2$&$.$&$.$&$.$&$-2$&$2$&$2$&$-2$&$-2$&$-2$\\
$4/2$&$32$&$.$&$.$&$-32$&$-32$&$32$&$-2$&$-2$&$.$&$.$&$.$&$2$&$2$&$2$&$-2$&$-2$&$2$\\
$5/2$&$48$&$.$&$.$&$48$&$-48$&$-48$&$.$&$.$&$.$&$.$&$.$&$.$&$.$&$.$&$.$&$.$&$.$\\
$6/2$&$64$&$.$&$.$&$-64$&$-64$&$64$&$2$&$2$&$.$&$.$&$.$&$-2$&$-2$&$-2$&$2$&$2$&$-2$\\
$7/2$&$96$&$.$&$.$&$96$&$-96$&$-96$&$.$&$.$&$.$&$.$&$.$&$.$&$.$&$.$&$.$&$.$&$.$\\
$8/2$&$132$&$.$&$.$&$-132$&$-132$&$132$&$.$&$.$&$.$&$4$&$.$&$.$&$.$&$.$&$.$&$.$&$.$\\
$9/2$&$184$&$.$&$.$&$184$&$-184$&$-184$&$2$&$2$&$.$&$.$&$.$&$-2$&$2$&$2$&$-2$&$-2$&$-2$\\
\hline
\end{tabular}}
\end{table}

\begin{table}[h!]\centering
\caption{Coefficients of $H_{(g,h),3}$.}
\resizebox{\columnwidth}{!}{
\begin{tabular}{r|rrrrrrrrrrrrrrrrr}
$$&$1a$&$2a$&$2b$&$2c$&$2d$&$2e$&$3a$&$3b$&$4a$&$4b$&$4c$&$6a$&$6b$&$6c$&$6d$&$6e$&$6f$\\
\hline$3/16$&$4$&$.$&$.$&$-4$&$4$&$-4$&$-1$&$-1$&$.$&$.$&$4$&$-1$&$1$&$1$&$1$&$1$&$-1$\\
$11/16$&$4$&$.$&$.$&$4$&$4$&$4$&$2$&$2$&$4$&$.$&$.$&$2$&$2$&$2$&$2$&$2$&$2$\\
$19/16$&$12$&$.$&$.$&$-12$&$12$&$-12$&$.$&$.$&$.$&$.$&$-4$&$.$&$.$&$.$&$.$&$.$&$.$\\
$27/16$&$16$&$.$&$.$&$16$&$16$&$16$&$-1$&$-1$&$.$&$.$&$.$&$-1$&$-1$&$-1$&$-1$&$-1$&$-1$\\
$35/16$&$24$&$.$&$.$&$-24$&$24$&$-24$&$.$&$.$&$.$&$.$&$8$&$.$&$.$&$.$&$.$&$.$&$.$\\
$43/16$&$36$&$.$&$.$&$36$&$36$&$36$&$.$&$.$&$4$&$.$&$.$&$.$&$.$&$.$&$.$&$.$&$.$\\
$51/16$&$56$&$.$&$.$&$-56$&$56$&$-56$&$-2$&$-2$&$.$&$.$&$-8$&$-2$&$2$&$2$&$2$&$2$&$-2$\\
$59/16$&$76$&$.$&$.$&$76$&$76$&$76$&$2$&$2$&$-4$&$.$&$.$&$2$&$2$&$2$&$2$&$2$&$2$\\
$67/16$&$108$&$.$&$.$&$-108$&$108$&$-108$&$.$&$.$&$.$&$.$&$12$&$.$&$.$&$.$&$.$&$.$&$.$\\
$75/16$&$148$&$.$&$.$&$148$&$148$&$148$&$-1$&$-1$&$4$&$.$&$.$&$-1$&$-1$&$-1$&$-1$&$-1$&$-1$\\
\hline
\end{tabular}}
\end{table}

\newpage

\begin{table}[h!]\centering
\caption{Decompositions for the first few $K^g_{1,\a}$.}
\resizebox{\columnwidth}{!}{
\begin{tabular}{r|rrrrrrrrrrrrrrr}
$$&$\chi_{1}$&$\chi_{2}$&$\chi_{3}$&$\chi_{4}$&$\chi_{5}$&$\chi_{6}$&$\chi_{7}$&$\chi_{8}$&$\chi_{9}$&$\chi_{10}$&$\chi_{11}$&$\chi_{12}$&$\chi_{13}$&$\chi_{14}$&$\chi_{15}$\\
\hline$3/16$&$.$&$.$&$.$&$1$&$1$&$.$&$.$&$.$&$.$&$.$&$.$&$.$&$.$&$.$&$.$\\
$11/16$&$.$&$2$&$.$&$.$&$.$&$.$&$.$&$.$&$.$&$.$&$.$&$.$&$.$&$.$&$.$\\
$19/16$&$.$&$.$&$.$&$.$&$.$&$.$&$.$&$.$&$.$&$.$&$.$&$.$&$2$&$.$&$.$\\
$27/16$&$.$&$.$&$.$&$.$&$.$&$1$&$1$&$.$&$.$&$.$&$.$&$.$&$.$&$2$&$.$\\
$35/16$&$2$&$.$&$.$&$2$&$2$&$.$&$.$&$.$&$.$&$.$&$.$&$.$&$2$&$.$&$.$\\
$43/16$&$.$&$2$&$.$&$.$&$.$&$2$&$2$&$.$&$.$&$.$&$.$&$.$&$.$&$4$&$.$\\
$51/16$&$.$&$.$&$.$&$2$&$2$&$.$&$.$&$.$&$.$&$.$&$.$&$.$&$8$&$.$&$.$\\
$59/16$&$.$&$4$&$.$&$.$&$.$&$2$&$2$&$.$&$.$&$.$&$.$&$.$&$.$&$10$&$.$\\
$67/16$&$6$&$.$&$.$&$6$&$6$&$.$&$.$&$.$&$.$&$.$&$.$&$.$&$12$&$.$&$.$\\
$75/16$&$.$&$6$&$.$&$.$&$.$&$7$&$7$&$.$&$.$&$.$&$.$&$.$&$.$&$18$&$.$\\
\hline
\end{tabular}}
\end{table}

\begin{table}[h!]\centering
\caption{Decompositions for the first few $K^g_{2,\a}$.}
\resizebox{\columnwidth}{!}{
\begin{tabular}{r|rrrrrrrrrrrrrrr}
$$&$\chi_{1}$&$\chi_{2}$&$\chi_{3}$&$\chi_{4}$&$\chi_{5}$&$\chi_{6}$&$\chi_{7}$&$\chi_{8}$&$\chi_{9}$&$\chi_{10}$&$\chi_{11}$&$\chi_{12}$&$\chi_{13}$&$\chi_{14}$&$\chi_{15}$\\
\hline$0/2$&$.$&$.$&$1$&$.$&$.$&$.$&$.$&$.$&$.$&$.$&$.$&$.$&$.$&$.$&$.$\\
$1/2$&$.$&$.$&$.$&$.$&$.$&$.$&$.$&$.$&$.$&$2$&$.$&$.$&$.$&$.$&$.$\\
$2/2$&$.$&$.$&$.$&$.$&$.$&$.$&$.$&$.$&$.$&$.$&$.$&$.$&$.$&$.$&$2$\\
$3/2$&$.$&$.$&$.$&$.$&$.$&$.$&$.$&$.$&$.$&$.$&$2$&$2$&$.$&$.$&$.$\\
$4/2$&$.$&$.$&$.$&$.$&$.$&$.$&$.$&$2$&$2$&$.$&$.$&$.$&$.$&$.$&$4$\\
$5/2$&$.$&$.$&$.$&$.$&$.$&$.$&$.$&$.$&$.$&$4$&$4$&$4$&$.$&$.$&$.$\\
$6/2$&$.$&$.$&$4$&$.$&$.$&$.$&$.$&$2$&$2$&$.$&$.$&$.$&$.$&$.$&$8$\\
$7/2$&$.$&$.$&$.$&$.$&$.$&$.$&$.$&$.$&$.$&$8$&$8$&$8$&$.$&$.$&$.$\\
$8/2$&$.$&$.$&$6$&$.$&$.$&$.$&$.$&$6$&$6$&$.$&$.$&$.$&$.$&$.$&$16$\\
$9/2$&$.$&$.$&$.$&$.$&$.$&$.$&$.$&$.$&$.$&$14$&$16$&$16$&$.$&$.$&$.$\\
\hline
\end{tabular}}
\end{table}

\begin{table}[h!]\centering
\caption{Decompositions for the first few $K^g_{3,\a}$.}
\resizebox{\columnwidth}{!}{
\begin{tabular}{r|rrrrrrrrrrrrrrr}
$$&$\chi_{1}$&$\chi_{2}$&$\chi_{3}$&$\chi_{4}$&$\chi_{5}$&$\chi_{6}$&$\chi_{7}$&$\chi_{8}$&$\chi_{9}$&$\chi_{10}$&$\chi_{11}$&$\chi_{12}$&$\chi_{13}$&$\chi_{14}$&$\chi_{15}$\\
\hline$3/16$&$.$&$.$&$.$&$.$&$.$&$1$&$1$&$.$&$.$&$.$&$.$&$.$&$.$&$.$&$.$\\
$11/16$&$2$&$.$&$.$&$.$&$.$&$.$&$.$&$.$&$.$&$.$&$.$&$.$&$.$&$.$&$.$\\
$19/16$&$.$&$.$&$.$&$.$&$.$&$.$&$.$&$.$&$.$&$.$&$.$&$.$&$.$&$2$&$.$\\
$27/16$&$.$&$.$&$.$&$1$&$1$&$.$&$.$&$.$&$.$&$.$&$.$&$.$&$2$&$.$&$.$\\
$35/16$&$.$&$2$&$.$&$.$&$.$&$2$&$2$&$.$&$.$&$.$&$.$&$.$&$.$&$2$&$.$\\
$43/16$&$2$&$.$&$.$&$2$&$2$&$.$&$.$&$.$&$.$&$.$&$.$&$.$&$4$&$.$&$.$\\
$51/16$&$.$&$.$&$.$&$.$&$.$&$2$&$2$&$.$&$.$&$.$&$.$&$.$&$.$&$8$&$.$\\
$59/16$&$4$&$.$&$.$&$2$&$2$&$.$&$.$&$.$&$.$&$.$&$.$&$.$&$10$&$.$&$.$\\
$67/16$&$.$&$6$&$.$&$.$&$.$&$6$&$6$&$.$&$.$&$.$&$.$&$.$&$.$&$12$&$.$\\
$75/16$&$6$&$.$&$.$&$7$&$7$&$.$&$.$&$.$&$.$&$.$&$.$&$.$&$18$&$.$&$.$\\
\hline
\end{tabular}}
\end{table}

\newpage

\newpage

\subsection[$\ell=6+3$, 3B twist]{$\boldsymbol{\ell=6+3}$, 3B twist}

\begin{table}[h!]\centering
\caption{Projective character table for $C_{G^{(6+3)}}(g)$.}
\begin{tabular}{r|rrrrrrrrr}
$$&$1a$&$2a$&$3a$&$3b$&$3c$&$3d$&$3e$&$6a$&$6b$\\
\hline$\chi_{1}$&$1$&$1$&$1$&$1$&$1$&$1$&$1$&$1$&$1$\\
$\chi_{2}$&$1$&$-1$&$1$&$1$&$1$&$1$&$1$&$-1$&$-1$\\
$\chi_{3}$&$1$&$-1$&$\zeta_{3}^{2}$&$\zeta_{3}$&$\zeta_{3}$&$1$&$\zeta_{3}^{2}$&$\zeta_{6}$&$\zeta_{6}^{5}$\\
$\chi_{4}$&$1$&$-1$&$\zeta_{3}$&$\zeta_{3}^{2}$&$\zeta_{3}^{2}$&$1$&$\zeta_{3}$&$\zeta_{6}^{5}$&$\zeta_{6}$\\
$\chi_{5}$&$1$&$1$&$\zeta_{3}^{2}$&$\zeta_{3}$&$\zeta_{3}$&$1$&$\zeta_{3}^{2}$&$\zeta_{3}^{2}$&$\zeta_{3}$\\
$\chi_{6}$&$1$&$1$&$\zeta_{3}$&$\zeta_{3}^{2}$&$\zeta_{3}^{2}$&$1$&$\zeta_{3}$&$\zeta_{3}$&$\zeta_{3}^{2}$\\
$\chi_{7}$&$2$&$.$&$2$&$2$&$-1$&$-1$&$-1$&$.$&$.$\\
$\chi_{8}$&$2$&$.$&$2\zeta_{3}$&$2\zeta_{3}^{2}$&$\zeta_{6}$&$-1$&$\zeta_{6}^{5}$&$.$&$.$\\
$\chi_{9}$&$2$&$.$&$2\zeta_{3}^{2}$&$2\zeta_{3}$&$\zeta_{6}^{5}$&$-1$&$\zeta_{6}$&$.$&$.$\\
\hline
\end{tabular}
\end{table}

\begin{table}[h!]\centering
\caption{Coefficients of $H_{(g,h),1}$.}
\begin{tabular}{r|rrrrrrrrr}
$$&$1a$&$2a$&$3a$&$3b$&$3c$&$3d$&$3e$&$6a$&$6b$\\
\hline$5/72$&$2$&$2$&$2$&$2$&$2$&$2$&$2$&$2$&$2$\\
$29/72$&$2$&$-2$&$2\zeta_{3}$&$2\zeta_{3}^{2}$&$2\zeta_{3}^{2}$&$2$&$2\zeta_{3}$&$2\zeta_{6}^{5}$&$2\zeta_{6}$\\
$53/72$&$4$&$.$&$4\zeta_{3}^{2}$&$4\zeta_{3}$&$4\zeta_{3}$&$4$&$4\zeta_{3}^{2}$&$.$&$.$\\
$77/72$&$2$&$-2$&$2$&$2$&$2$&$2$&$2$&$-2$&$-2$\\
$101/72$&$6$&$2$&$6\zeta_{3}$&$6\zeta_{3}^{2}$&$6\zeta_{3}^{2}$&$6$&$6\zeta_{3}$&$2\zeta_{3}$&$2\zeta_{3}^{2}$\\
$125/72$&$6$&$-2$&$6\zeta_{3}^{2}$&$6\zeta_{3}$&$6\zeta_{3}$&$6$&$6\zeta_{3}^{2}$&$2\zeta_{6}$&$2\zeta_{6}^{5}$\\
$149/72$&$8$&$.$&$8$&$8$&$8$&$8$&$8$&$.$&$.$\\
$173/72$&$8$&$.$&$8\zeta_{3}$&$8\zeta_{3}^{2}$&$8\zeta_{3}^{2}$&$8$&$8\zeta_{3}$&$.$&$.$\\
$197/72$&$12$&$4$&$12\zeta_{3}^{2}$&$12\zeta_{3}$&$12\zeta_{3}$&$12$&$12\zeta_{3}^{2}$&$4\zeta_{3}^{2}$&$4\zeta_{3}$\\
$221/72$&$14$&$-2$&$14$&$14$&$14$&$14$&$14$&$-2$&$-2$\\
\hline
\end{tabular}
\end{table}

\begin{table}[h!]\centering
\caption{Coefficients of $H_{(g,h),3}$.}
\begin{tabular}{r|rrrrrrrrr}
$$&$1a$&$2a$&$3a$&$3b$&$3c$&$3d$&$3e$&$6a$&$6b$\\
\hline$5/72$&$4$&$.$&$4\zeta_{3}$&$4\zeta_{3}^{2}$&$2\zeta_{6}$&$-2$&$2\zeta_{6}^{5}$&$.$&$.$\\
$29/72$&$4$&$.$&$4\zeta_{3}^{2}$&$4\zeta_{3}$&$2\zeta_{6}^{5}$&$-2$&$2\zeta_{6}$&$.$&$.$\\
$53/72$&$8$&$.$&$8$&$8$&$-4$&$-4$&$-4$&$.$&$.$\\
$77/72$&$4$&$.$&$4\zeta_{3}$&$4\zeta_{3}^{2}$&$2\zeta_{6}$&$-2$&$2\zeta_{6}^{5}$&$.$&$.$\\
$101/72$&$12$&$.$&$12\zeta_{3}^{2}$&$12\zeta_{3}$&$6\zeta_{6}^{5}$&$-6$&$6\zeta_{6}$&$.$&$.$\\
$125/72$&$12$&$.$&$12$&$12$&$-6$&$-6$&$-6$&$.$&$.$\\
$149/72$&$16$&$.$&$16\zeta_{3}$&$16\zeta_{3}^{2}$&$8\zeta_{6}$&$-8$&$8\zeta_{6}^{5}$&$.$&$.$\\
$173/72$&$16$&$.$&$16\zeta_{3}^{2}$&$16\zeta_{3}$&$8\zeta_{6}^{5}$&$-8$&$8\zeta_{6}$&$.$&$.$\\
$197/72$&$24$&$.$&$24$&$24$&$-12$&$-12$&$-12$&$.$&$.$\\
$221/72$&$28$&$.$&$28\zeta_{3}$&$28\zeta_{3}^{2}$&$14\zeta_{6}$&$-14$&$14\zeta_{6}^{5}$&$.$&$.$\\
\hline
\end{tabular}
\end{table}

\begin{table}[h!]\centering
\caption{Decompositions for the first few $K^g_{1,\a}$.}
\begin{tabular}{r|rrrrrrrrr}
$$&$\chi_{1}$&$\chi_{2}$&$\chi_{3}$&$\chi_{4}$&$\chi_{5}$&$\chi_{6}$&$\chi_{7}$&$\chi_{8}$&$\chi_{9}$\\
\hline$5/72$&$2$&$.$&$.$&$.$&$.$&$.$&$.$&$.$&$.$\\
$29/72$&$.$&$.$&$.$&$2$&$.$&$.$&$.$&$.$&$.$\\
$53/72$&$.$&$.$&$2$&$.$&$2$&$.$&$.$&$.$&$.$\\
$77/72$&$.$&$2$&$.$&$.$&$.$&$.$&$.$&$.$&$.$\\
$101/72$&$.$&$.$&$.$&$2$&$.$&$4$&$.$&$.$&$.$\\
$125/72$&$.$&$.$&$4$&$.$&$2$&$.$&$.$&$.$&$.$\\
$149/72$&$4$&$4$&$.$&$.$&$.$&$.$&$.$&$.$&$.$\\
$173/72$&$.$&$.$&$.$&$4$&$.$&$4$&$.$&$.$&$.$\\
$197/72$&$.$&$.$&$4$&$.$&$8$&$.$&$.$&$.$&$.$\\
$221/72$&$6$&$8$&$.$&$.$&$.$&$.$&$.$&$.$&$.$\\
\hline
\end{tabular}
\end{table}

\begin{table}[h!]\centering
\caption{Decompositions for the first few $K^g_{3,\a}$.}
\begin{tabular}{r|rrrrrrrrr}
$$&$\chi_{1}$&$\chi_{2}$&$\chi_{3}$&$\chi_{4}$&$\chi_{5}$&$\chi_{6}$&$\chi_{7}$&$\chi_{8}$&$\chi_{9}$\\
\hline$5/72$&$.$&$.$&$.$&$.$&$.$&$.$&$.$&$2$&$.$\\
$29/72$&$.$&$.$&$.$&$.$&$.$&$.$&$.$&$.$&$2$\\
$53/72$&$.$&$.$&$.$&$.$&$.$&$.$&$4$&$.$&$.$\\
$77/72$&$.$&$.$&$.$&$.$&$.$&$.$&$.$&$2$&$.$\\
$101/72$&$.$&$.$&$.$&$.$&$.$&$.$&$.$&$.$&$6$\\
$125/72$&$.$&$.$&$.$&$.$&$.$&$.$&$6$&$.$&$.$\\
$149/72$&$.$&$.$&$.$&$.$&$.$&$.$&$.$&$8$&$.$\\
$173/72$&$.$&$.$&$.$&$.$&$.$&$.$&$.$&$.$&$8$\\
$197/72$&$.$&$.$&$.$&$.$&$.$&$.$&$12$&$.$&$.$\\
$221/72$&$.$&$.$&$.$&$.$&$.$&$.$&$.$&$14$&$.$\\
\hline
\end{tabular}
\end{table}

\newpage

\subsection*{Acknowledgements}

We are grateful to John Duncan, Simon Lentner, Terry Gannon, Erik Verlinde for helpful discussions. We especially John Duncan for many of the group descriptions and for helpful comments on an earlier version of the manuscript. The work of M.C.\ and D.W.\ was supported by ERC starting grant H2020 ERC StG 2014.

\pdfbookmark[1]{References}{ref}
\LastPageEnding


\begin{thebibliography}{99}
\footnotesize\itemsep=0pt

\bibitem{Altschuler:2003nn}
Altschuler D., Coste A., Maillard J.M., Representation theory of twisted group
 double, \textit{Ann. Fond. Louis de Broglie} \textbf{29} (2004), 681--694,
 \href{https://arxiv.org/abs/hep-th/0309257}{arXiv:hep-th/0309257}.

\bibitem{Anagiannis:2017src}
Anagiannis V., Cheng M.C.N., Harrison S.M., $K3$ elliptic genus and an umbral
 moonshine module, \href{https://arxiv.org/abs/1709.01952}{arXiv:1709.01952}.

\bibitem{Bantay:1990tr}
B\'{a}ntay P., Orbifolds, {H}opf algebras, and the moonshine, \href{https://doi.org/10.1007/BF00403544}{\textit{Lett.
 Math. Phys.}} \textbf{22} (1991), 187--194.

\bibitem{Bor_PNAS}
Borcherds R.E., Vertex algebras, {K}ac--{M}oody algebras, and the monster,
 \href{https://doi.org/10.1073/pnas.83.10.3068}{\textit{Proc. Nat. Acad. Sci. USA}} \textbf{83} (1986), 3068--3071.

\bibitem{borcherds_monstrous}
Borcherds R.E., Monstrous moonshine and monstrous {L}ie superalgebras,
 \href{https://doi.org/10.1007/BF01232032}{\textit{Invent. Math.}} \textbf{109} (1992), 405--444.

\bibitem{MR1324339}
Brown K.S., Cohomology of groups, \textit{Graduate Texts in Mathematics},
 Vol.~87, \href{https://doi.org/10.1007/978-1-4684-9327-6}{Springer-Verlag}, New York, 1994.

\bibitem{Cappelli:1987xt}
Cappelli A., Itzykson C., Zuber J.B., The {${\rm A}$}-{${\rm D}$}-{${\rm E}$}
 classification of minimal and {$A^{(1)}_1$} conformal invariant theories,
 \href{https://doi.org/10.1007/BF01221394}{\textit{Comm. Math. Phys.}} \textbf{113} (1987), 1--26.

\bibitem{Carnahan2008}
Carnahan S., Generalized moonshine {I}: {G}enus-zero functions, \href{https://doi.org/10.2140/ant.2010.4.649}{\textit{Algebra
 Number Theory}} \textbf{4} (2010), 649--679, \href{https://arxiv.org/abs/0812.3440}{arXiv:0812.3440}.

\bibitem{MR2904095}
Carnahan S., Generalized moonshine, {II}:~{B}orcherds products, \href{https://doi.org/10.1215/00127094-1548416}{\textit{Duke
 Math.~J.}} \textbf{161} (2012), 893--950, \href{https://arxiv.org/abs/0908.4223}{arXiv:0908.4223}.

\bibitem{carnahan12}
Carnahan S., Generalized moonshine, {IV}: {M}onstrous {L}ie algebras,
 \href{https://arxiv.org/abs/1208.6254}{arXiv:1208.6254}.

\bibitem{carnahan_miyamoto}
Carnahan S., Miyamoto M., Regularity of fixed-point vertex operator
 subalgebras, \href{https://arxiv.org/abs/1603.05645}{arXiv:1603.05645}.

\bibitem{Cheng2010_1}
Cheng M.C.N., {$K3$} surfaces, {${\mathcal N}=4$} dyons and the {M}athieu group
 {$M_{24}$}, \href{https://doi.org/10.4310/CNTP.2010.v4.n4.a2}{\textit{Commun. Number Theory Phys.}} \textbf{4} (2010), 623--657,
 \href{https://arxiv.org/abs/1005.5415}{arXiv:1005.5415}.

\bibitem{Cheng2011}
Cheng M.C.N., Duncan J.F.R., On {R}ademacher sums, the largest {M}athieu group
 and the holographic modularity of moonshine, \href{https://doi.org/10.4310/CNTP.2012.v6.n3.a4}{\textit{Commun. Number Theory
 Phys.}} \textbf{6} (2012), 697--758, \href{https://arxiv.org/abs/1110.3859}{arXiv:1110.3859}.

\bibitem{2012arXiv1212.0906C}
Cheng M.C.N., Duncan J.F.R., On the discrete groups of {M}athieu moonshine, in
 Perspectives in Representation Theory, \textit{Contemp. Math.}, Vol.~610,
 \href{https://doi.org/10.1090/conm/610/12195}{Amer. Math. Soc.}, Providence, RI, 2014, 65--78, \href{https://arxiv.org/abs/1212.0906}{arXiv:1212.0906}.

\bibitem{omjt}
Cheng M.C.N., Duncan J.F.R., Optimal mock Jacobi theta functions,
 \href{https://arxiv.org/abs/1605.04480}{arXiv:1605.04480}.


\bibitem{Cheng:2017grj}
Cheng M.C.N., Duncan J.F.R., Meromorphic {J}acobi forms of half-integral index
 and umbral moonshine modules, \href{https://arxiv.org/abs/1707.01336}{arXiv:1707.01336}.

\bibitem{UM}
Cheng M.C.N., Duncan J.F.R., Harvey J.A., Umbral moonshine, \href{https://doi.org/10.4310/CNTP.2014.v8.n2.a1}{\textit{Commun.
 Number Theory Phys.}} \textbf{8} (2014), 101--242, \href{https://arxiv.org/abs/1204.2779}{arXiv:1204.2779}.

\bibitem{MUM}
Cheng M.C.N., Duncan J.F.R., Harvey J.A., Umbral moonshine and the {N}iemeier
 lattices, \href{https://doi.org/10.1186/2197-9847-1-3}{\textit{Res. Math. Sci.}} \textbf{1} (2014), Art.~3, 81~pages,
 \href{https://arxiv.org/abs/1307.5793}{arXiv:1307.5793}.

\bibitem{mumcorr}
Cheng M.C.N., Duncan J.F.R., Harvey J.A., Weight one {J}acobi forms and umbral
 moonshine, \href{https://doi.org/10.1088/1751-8121/aaa819}{\textit{J.~Phys.~A: Math. Theor.}} \textbf{51} (2018), 104002,
 37~pages, \href{https://arxiv.org/abs/1703.03968}{arXiv:1703.03968}.

\bibitem{UMLG}
Cheng M.C.N., Ferrari F., Harrison S.M., Paquette N.M., Landau--{G}inzburg
 orbifolds and symmetries of ${K}3$ {CFT}s, \href{https://doi.org/10.1007/JHEP01(2017)046}{\textit{J.~High Energy Phys.}}
 \textbf{2017} (2017), no.~1, 046, 49~pages, \href{https://arxiv.org/abs/1512.04942}{arXiv:1512.04942}.

\bibitem{UMk3}
Cheng M.C.N., Harrison S., Umbral moonshine and $K3$ surfaces, \href{https://doi.org/10.1007/s00220-015-2398-5}{\textit{Comm.
 Math. Phys.}} \textbf{339} (2015), 221--261, \href{https://arxiv.org/abs/1406.0619}{arXiv:1406.0619}.

\bibitem{conway_norton}
Conway J.H., Norton S.P., Monstrous moonshine, \href{https://doi.org/10.1112/blms/11.3.308}{\textit{Bull. London Math. Soc.}}
 \textbf{11} (1979), 308--339.

\bibitem{Dabholkar:2012nd}
Dabholkar A., Murthy S., Zagier D., Quantum black holes, wall crossing, and
 mock modular forms, \href{https://arxiv.org/abs/1208.4074}{arXiv:1208.4074}.

\bibitem{Roche:1990hs}
Dijkgraaf R., Pasquier V., Roche P., Quasi {H}opf algebras, group cohomology
 and orbifold models, \href{https://doi.org/10.1016/0920-5632(91)90123-V}{\textit{Nuclear Phys.~B Proc. Suppl.}} \textbf{18B}
 (1990), 60--72.

\bibitem{Dijkgraaf:1989hb}
Dijkgraaf R., Vafa C., Verlinde E., Verlinde H., The operator algebra of
 orbifold models, \href{https://doi.org/10.1007/BF01238812}{\textit{Comm. Math. Phys.}} \textbf{123} (1989), 485--526.

\bibitem{Dijkgraaf:1989pz}
Dijkgraaf R., Witten E., Topological gauge theories and group cohomology,
 \href{https://doi.org/10.1007/BF02096988}{\textit{Comm. Math. Phys.}} \textbf{129} (1990), 393--429.

\bibitem{Dixon:1988qd}
Dixon L., Ginsparg P., Harvey J., Beauty and the beast: superconformal symmetry
 in a {M}onster module, \href{https://doi.org/10.1007/BF01217740}{\textit{Comm. Math. Phys.}} \textbf{119} (1988),
 221--241.

\bibitem{Dixon:1985jw}
Dixon L., Harvey J.A., Vafa C., Witten E., Strings on orbifolds,
 \href{https://doi.org/10.1016/0550-3213(85)90593-0}{\textit{Nuclear Phys.~B}} \textbf{261} (1985), 678--686.

\bibitem{Dixon:1986jc}
Dixon L., Harvey J.A., Vafa C., Witten E., Strings on orbifolds.~{II},
 \href{https://doi.org/10.1016/0550-3213(86)90287-7}{\textit{Nuclear Phys.~B}} \textbf{274} (1986), 285--314.

\bibitem{Dong2000}
Dong C., Li H., Mason G., Modular-invariance of trace functions in orbifold
 theory and generalized moonshine, \href{https://doi.org/10.1007/s002200000242}{\textit{Comm. Math. Phys.}} \textbf{214}
 (2000), 1--56, \href{https://arxiv.org/abs/q-alg/9703016}{arXiv:q-alg/9703016}.

\bibitem{MR934283}
Drinfel'd V.G., Quantum groups, in Proceedings of the {I}nternational
 {C}ongress of {M}athematicians, {V}ols.~1,~2 ({B}erkeley, {C}alif., 1986),
 Amer. Math. Soc., Providence, RI, 1987, 798--820.

\bibitem{proof}
Duncan J.F.R., Griffin M.J., Ono K., Proof of the umbral moonshine conjecture,
 \href{https://doi.org/10.1186/s40687-015-0044-7}{\textit{Res. Math. Sci.}} \textbf{2} (2015), Art.~26, 47~pages,
 \href{https://arxiv.org/abs/1406.0571}{arXiv:1406.0571}.

\bibitem{Duncan:2014tya}
Duncan J.F.R., Harvey J.A., The umbral moonshine module for the unique
 unimodular {N}iemeier root system, \href{https://doi.org/10.2140/ant.2017.11.505}{\textit{Algebra Number Theory}} \textbf{11}
 (2017), 505--535, \href{https://arxiv.org/abs/1412.8191}{arXiv:1412.8191}.

\bibitem{Duncan:2017bhh}
Duncan J.F.R., O'Desky A., Super vertex algebras, meromorphic {J}acobi forms
 and umbral moonshine, \href{https://doi.org/10.1016/j.jalgebra.2018.08.017}{\textit{J.~Algebra}} \textbf{515} (2018), 389--407,
 \href{https://arxiv.org/abs/1705.09333}{arXiv:1705.09333}.

\bibitem{Eguchi2010a}
Eguchi T., Hikami K., Note on twisted elliptic genus of {$K3$} surface,
 \href{https://doi.org/10.1016/j.physletb.2010.10.017}{\textit{Phys. Lett.~B}} \textbf{694} (2011), 446--455, \href{https://arxiv.org/abs/1008.4924}{arXiv:1008.4924}.

\bibitem{Eguchi2010}
Eguchi T., Ooguri H., Tachikawa Y., Notes on the {$K3$} surface and the
 {M}athieu group {$M_{24}$}, \href{https://doi.org/10.1080/10586458.2011.544585}{\textit{Exp. Math.}} \textbf{20} (2011), 91--96,
 \href{https://arxiv.org/abs/1004.0956}{arXiv:1004.0956}.

\bibitem{eichler_zagier}
Eichler M., Zagier D., The theory of {J}acobi forms, \textit{Progress in
 Mathematics}, Vol.~55, \href{https://doi.org/10.1007/978-1-4684-9162-3}{Birkh\"{a}user Boston}, Inc., Boston, MA, 1985.

\bibitem{FLM}
Frenkel I., Lepowsky J., Meurman A., Vertex operator algebras and the monster,
 \textit{Pure and Applied Mathematics}, Vol.~134, Academic Press, Inc.,
 Boston, MA, 1988.

\bibitem{FLMPNAS}
Frenkel I.B., Lepowsky J., Meurman A., A natural representation of the
 {F}ischer--{G}riess monster with the modular function {$J$} as character,
 \href{https://doi.org/10.1073/pnas.81.10.3256}{\textit{Proc. Nat. Acad. Sci. USA}} \textbf{81} (1984), 3256--3260.

\bibitem{FLMBerk}
Frenkel I.B., Lepowsky J., Meurman A., A moonshine module for the monster, in
 Vertex Operators in Mathematics and Physics ({B}erkeley, {C}alif., 1983),
 \textit{Math. Sci. Res. Inst. Publ.}, Vol.~3, \href{https://doi.org/10.1007/978-1-4613-9550-8_12}{Springer}, New York, 1985,
 231--273.

\bibitem{Gaberdiel2010}
Gaberdiel M.R., Hohenegger S., Volpato R., Mathieu twining characters for
 {$K3$}, \href{https://doi.org/10.1007/JHEP09(2010)058}{\textit{J.~High Energy Phys.}} \textbf{2010} (2010), no.~9, 058,
 20~pages, \href{https://arxiv.org/abs/1006.0221}{arXiv:1006.0221}.

\bibitem{Gaberdiel2010a}
Gaberdiel M.R., Hohenegger S., Volpato R., Mathieu moonshine in the elliptic
 genus of {$K3$}, \href{https://doi.org/10.1007/JHEP10(2010)062}{\textit{J.~High Energy Phys.}} \textbf{2010} (2010), no.~10,
 062, 24~pages, \href{https://arxiv.org/abs/1008.3778}{arXiv:1008.3778}.

\bibitem{Gaberdiel:2012gf}
Gaberdiel M.R., Persson D., Ronellenfitsch H., Volpato R., Generalized
 {M}athieu moonshine, \href{https://doi.org/10.4310/CNTP.2013.v7.n1.a5}{\textit{Commun. Number Theory Phys.}} \textbf{7} (2013),
 145--223, \href{https://arxiv.org/abs/1211.7074}{arXiv:1211.7074}.

\bibitem{Gannon:2012ck}
Gannon T., Much ado about {M}athieu, \href{https://doi.org/10.1016/j.aim.2016.06.014}{\textit{Adv. Math.}} \textbf{301} (2016),
 322--358, \href{https://arxiv.org/abs/1211.5531}{arXiv:1211.5531}.

\bibitem{MR3582425}
Griffin M.J., Mertens M.H., A proof of the {T}hompson moonshine conjecture,
 \href{https://doi.org/10.1186/s40687-016-0084-7}{\textit{Res. Math. Sci.}} \textbf{3} (2016), Art.~36, 32~pages,
 \href{https://arxiv.org/abs/1607.03078}{arXiv:1607.03078}.

\bibitem{Harvey:2013mda}
Harvey J.A., Murthy S., Moonshine in fivebrane spacetimes, \href{https://doi.org/10.1007/JHEP01(2014)146}{\textit{J.~High
 Energy Phys.}} \textbf{2014} (2014), no.~1, 146, 29~pages, \href{https://arxiv.org/abs/1307.7717}{arXiv:1307.7717}.

\bibitem{Harvey:2014cva}
Harvey J.A., Murthy S., Nazaroglu C., ADE double scaled little string theories,
 mock modular forms and umbral moonshine, \href{https://doi.org/10.1007/JHEP05(2015)126}{\textit{J.~High Energy Phys.}}
 \textbf{2015} (2015), no.~5, 126, 57~pages, \href{https://arxiv.org/abs/1410.6174}{arXiv:1410.6174}.

\bibitem{Hoehn03}
H\"ohn G., Generalized moonshine for the baby monster, 2003.

\bibitem{Kachru:2016ttg}
Kachru S., Paquette N.M., Volpato R., 3{D} string theory and umbral moonshine,
 \href{https://doi.org/10.1088/1751-8121/aa6e07}{\textit{J.~Phys.~A: Math. Theor.}} \textbf{50} (2017), 404003, 21~pages,
 \href{https://arxiv.org/abs/1603.07330}{arXiv:1603.07330}.

\bibitem{generalized_moonshine}
Mason G., Finite groups and modular functions (with an appendix by
 S.P.~Norton), in The {A}rcata {C}onference on {R}epresentations of {F}inite
 {G}roups ({A}rcata, {C}alif., 1986), \textit{Proc. Sympos. Pure Math.},
 Vol.~47, Amer. Math. Soc., Providence, RI, 1987, 181--210.

\bibitem{generalized_moonshine2}
Norton S., From moonshine to the monster, in Proceedings on {M}oonshine and
 Related Topics ({M}ontr\'{e}al, {QC}, 1999), \textit{CRM Proc. Lecture
 Notes}, Vol.~30, Amer. Math. Soc., Providence, RI, 2001, 163--171.

\bibitem{Paquette:2016xoo}
Paquette N.M., Persson D., Volpato R., Monstrous {BPS}-algebras and the
 superstring origin of moonshine, \href{https://doi.org/10.4310/CNTP.2016.v10.n3.a2}{\textit{Commun. Number Theory Phys.}}
 \textbf{10} (2016), 433--526, \href{https://arxiv.org/abs/1601.05412}{arXiv:1601.05412}.

\bibitem{MR0472707}
Shimura G., On modular forms of half integral weight, \href{https://doi.org/10.2307/1970831}{\textit{Ann. of Math.}}
 \textbf{97} (1973), 440--481.

\bibitem{MR1372731}
Tuite M.P., Generalised moonshine and abelian orbifold constructions, in
 Moonshine, the {M}onster, and Related Topics ({S}outh {H}adley, {MA}, 1994),
 \textit{Contemp. Math.}, Vol.~193, \href{https://doi.org/10.1090/conm/193/02380}{Amer. Math. Soc.}, Providence, RI, 1996,
 353--368, \href{https://arxiv.org/abs/hep-th/9412036}{arXiv:hep-th/9412036}.

\bibitem{Rad_whalen}
Whalen D., Vector-valued {R}ademacher sums and automorphic integrals,
 \href{https://arxiv.org/abs/1406.0571}{arXiv:1406.0571}.

\bibitem{Zhu_ModInv}
Zhu Y., Modular invariance of characters of vertex operator algebras,
 \href{https://doi.org/10.1090/S0894-0347-96-00182-8}{\textit{J.~Amer. Math. Soc.}} \textbf{9} (1996), 237--302.

\bibitem{Zwegers2008}
Zwegers S., Mock theta functions, Ph.D.~Thesis, {U}trecht {U}niversity, 2002,
 \href{https://arxiv.org/abs/0807.4834}{arXiv:0807.4834}.

\end{thebibliography}
\end{document}